\documentclass{article}
\usepackage[english]{babel}
\usepackage{amsmath,amssymb,graphicx,xcolor,latexsym}

\newcommand{\Omicron}{\mathrm{O}}
\newcommand{\infixand}{\text{ and }}
\newcommand{\mathd}{\mathrm{d}}

\newcommand{\tmmathbf}[1]{\ensuremath{\boldsymbol{#1}}}
\newcommand{\tmop}[1]{\ensuremath{\operatorname{#1}}}
\newcommand{\tmscript}[1]{\text{\scriptsize{$#1$}}}

\newtheorem{lemma}{Lemma}
\newtheorem{theorem}{Theorem}

\begin{document}

\title{The optimal bound $e$ in a dyadic version of Uchiyama's Lemma}

\author{J. Fladung and S. Petermichl}

\maketitle

\begin{abstract}
  The best known constant in Uchiyama's Lemma is $e$. A conjecture states that
  this cannot be improved. We show that the constant $e$ also stands in a
  dyadic version of Uchiyama's Lemma. Further, we prove that in the dyadic
  case, the constant $e$ can indeed not be improved. We deduce a dyadic version of the reproducing kernel thesis for the embedding theorem.
\end{abstract}

\section{Context of our subject}

Carleson embedding theorems are at the heart of several branches of analysis. In complex analysis, one finds them stated for analytic maps in the unit disk with an embedding of the analytic Hardy space $H^2$. In other parts of analysis one finds $L^2$ versions with harmonic functions or their generalizations, with extensive meaning and applications to PDE. Since classical embeddings are concerned with square integrals and are often formulated with packing conditions on the measure, there is no difference if one considers analytic or harmonic maps for the embedding itself - that is - as long as one does not care about sharp constants in the embedding.

The classical packing condition itself in the Carleson embedding theorem is not well adapted to considerations with sharp estimates, as the condition is purely geometric and thus outside the scope of the structure that governs the embedding itself. When formulating the packing condition by means of Brownian trajectories or Green potentials, meaningful constants can be obtained. This type of theorem is Uchiyama's Lemma in complex analysis and a formulation with Green potentials for the harmonic case.

The testing form on the other hand relies on testing the embedding on a relevant class of the desired embedding and is as such always structurally sound. In the complex formulation, for example, these functions are the reproducing kernels of the $H^2$ space.

Notably dyadic versions of Carleson embeddings have played a large role in recent years \cite{NTV}, \cite{NPTV}, \cite{CT} as they tend to contain the essential information.

In the classic dyadic form of the embedding theorem, these two different formulations (testing and packing) are almost indistinguishable. The test functions are characteristic functions of dyadic intervals - a renormalization of the reproducing kernel for dyadic extensions. One comes back to the packing condition trivially by testing on characteristic functions.

In the dyadic case the best constant for the Carleson embedding theorem with dyadic packing condition is 4 and it is known that this constant is sharp. The constant 4 is then inherited by the harmonic version in its formulation with the Green potential.

The best known constant in Uchiyama's Lemma is $e$. It is smaller than 4 thanks to the restriction to analytic maps and it has been conjectured that this constant is sharp  \cite{PTW}. In this paper we present a thoughtful dyadic model for the complex case and show that the embedding in this case still has a constant $e$ and - most importantly - that $e$ is sharp.

From our dyadic version of the analytic embedding we then proceed to prove the reproducing kernel thesis for embedding theorems in this setting by testing on reproducing kernels with a constant $3e$.

\paragraph{Uchiyama's Lemma and its dyadic version.}

The classical form is the following. A printed version can be found in Nikolskii's book \cite{N}.

\begin{lemma}[Uchiyama]\label{lemmauchiyama}
Let $\varphi$ be a non-positive subharmonic function in $\mathbb{D}$. Then the
 measure
\[ \mathd \nu = e^{\varphi} \Delta \varphi \log \frac{1}{| z |} \mathd A (z)
\]
is Carleson and moreover for $f \in H^2 (\mathbb{D})$ there holds the
embedding
\[ \int_{\mathbb{D}} | f (z) |^2 \mathd \nu (z) \leqslant \| f \|_{H^2}^2 . \]
\end{lemma}

A simple consequence, by estimating the exponential function from below and
scaling, is the following:

\begin{theorem}[Uchiyama]
  \label{corollaryuchiyama}Let $\varphi$ be a bounded non-positive subharmonic
  function on $\mathbb{D}$. Then for the measure
  \[ \mathd \mu (z) = \Delta \varphi (z) \log \frac{1}{| z |} \mathd A (z), \]
  $H^2 (\mathbb{D})$ embeds in $L^2 (\mu)$ with constant $e$: for all $f\in H^2 (\mathbb{D})$ there holds
  \[ \int_{\mathbb{D}} | f (z) |^2 \mathd \mu (z) \leqslant e \| \varphi
     \|_{\infty} \| f \|_{H^2}^2 . \]
\end{theorem}

Theorem \ref{corollaryuchiyama} is an embedding theorem, but it is still often refererred to as Uchiyama's Lemma.

\

The auxiliary function $(\varphi, f) \mapsto e^{\varphi} | f |^2$ implies the
Lemma \ref{lemmauchiyama} immediately through a use of Green's formula after computing
\[ \Delta (e^{\varphi} | f |^2) = e^{\varphi} \Delta \varphi | f |^2 + 4
   e^{\varphi} | \partial \varphi f + \partial f |^2 \geqslant e^{\varphi}
   \Delta \varphi | f |^2 . \]
The constant $e$ in Theorem \ref{corollaryuchiyama} was conjectured to be
optimal in \cite{PTW}, even though someone might suspect it an artifact of the proof: it
arises from the use of the function $(\varphi, f) \mapsto e^{\varphi} | f |^2$
and $e^{-1}=\min \{ e^{\varphi} : - 1 \leqslant \varphi \leqslant 0
\}$ if one assumes $\| \varphi \|_{\infty} = 1$. (The general case is then
obtained by scaling.)

\

As the subject of our note will show, the constant $e$ that arises can at least in some sense not be improved, as we consider a dyadic model for the continuous problem. In this formulation the constant $e$ still stands and cannot be improved.

\

To study a dyadic version of Uchiyama's lemma, we
define a `dyadic analyticity' via `analytic' (or orthognal) increments on the real part $u$ and the
imaginary part $v$ indexed by 4-adic intervals $I$, as illustrated by the picture below. It shows the typical
rotation one expects from pairs of conjugate functions so that we have the
`dyadic Cauchy Riemann equations'
\begin{equation}
\label{dyadicCR}
 \Delta^x u_I = \Delta^y v_I, \Delta^y u_I = - \Delta^x v_I.
\end{equation}

\raisebox{-0.5\height}{\includegraphics[width=11.1cm,height=6.75cm]{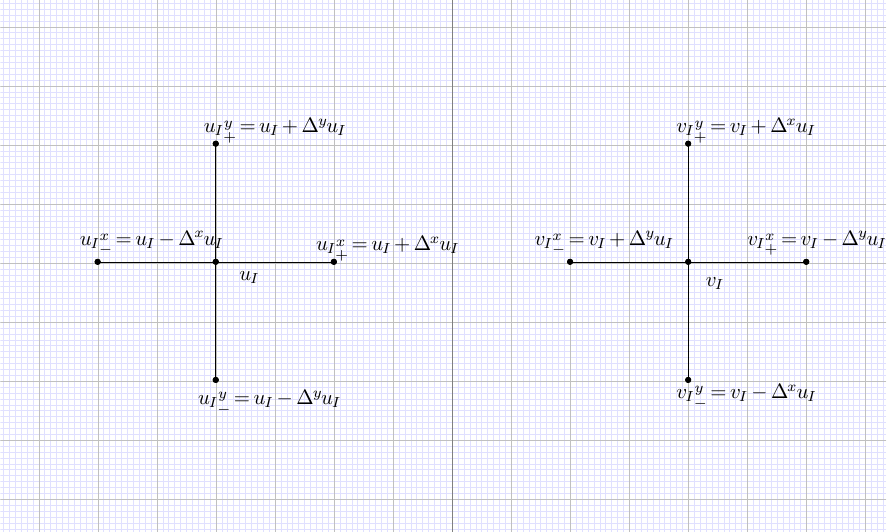}}

\

For both $u$ and $v$ this gives a 4-adic tree with the constraint that the averages of the left
two children equal the entry of the root and the same for the right two
children - with increments in their roles as rescaled $x$ and $y$ derivatives. Such a tree can be
identified with a sliced dyadic martingale, that is, a dyadic martingale with zero increments in every
second step. From a probabilistic
standpoint, $u$ and $v$ are thus sliced (0 increments in every second step)
orthogonal  (they have orthogonal increments) $L^2$ integrable (and therefore closable) dyadic martingales. Generally
we call $f = \tmop{Re} f + i \tmop{Im} f$ `dyadic analytic' if $f$ is a sliced
dyadic martingale whose real and imaginary parts also have dyadic Cauchy Riemann \eqref{dyadicCR}. See section \eqref{definitionsection} for the precise definitions.

The term `dyadic analytic' is a bit abusive, as we have a lot more flexibility with orthogonal martingales than with analytic maps. However, such models have proven precise examples of the Hilbert transform in the past, and have as such a natural link to analytic functions. As we will see below, the conjugate $v$ can be obtained from $u$ through the use of the operator $\mathcal{S}_0$ that acts directly on the Haar base through a rotation of neighbouring increments.
The operator $\mathcal{S}_0$ has been used as a
dyadic model for the Hilbert transform and was a decisive tool in \cite{DP}. Other dyadic operators have also served as dyadic form of the Hilbert transform, for example in \cite{P}, \cite{DPTV}.

\

In our model the dyadic Hardy space $H_{\tmop{dy}}^2$ consists of pairs $(u,
v)$ (or dyadic analytic $f = u + i v$) of sliced $L^2$ integrable dyadic
martingales with dyadic Cauchy-Riemann relations. It is a reproducing kernel Hilbert space with its inner product inherited from the inner product for complex valued $L^2$ functions.

\

We also have a sliced supermartingale $M=(M_I)_{I\in \mathcal{D}^4}$. 
Here we require
\[ M_I \geqslant \frac{M_{I_+^x} + M_{I_-^x}}{2} = \frac{M_{I_+^y} + M_{I_-^y}}{2}, \]
so that in particular there exists $\mu_I\geqslant 0$ so that for $\star = x, y$
\[ -\Delta^{\tmop{dy}} M_I:=- \Delta_{\star}^{\tmop{dy}} M_I := M_I - \frac{M_{I_+^{\star}} + M_{I_-^{\star}}}{2} =\frac{\mu_I}{|I|}. \]
If the $\mu_I = 0$ then we have again a sliced martingale. Depending on the sign one adapts, one can of course formulate everything using a non-positive sliced submartingale $M$.

\

The operator $\mathcal{S}_0$ here only acts on sliced dyadic martingales. While the slicing is well motivated so as to have $x$ and $y$ derivatives `attach in the same point', it is also essential if one is trying to use the Bellman function we used: already dropping the slicing condition on the supermartingale $M$ cannot be handled by the Bellman function we use in our proofs, even if only applied to sliced $u$ and $v$. See the remarks near the end. The main interest of this project was for us the proof of the optimality in the case of a suitable model for which the upper estimate still holds. Getting a blow up in a more restrictive setting for which an upper estimate can still be proved, is in some sense a stronger result.

\

The dyadic form of Uchiyama's lemma is this.

\begin{lemma}\label{dyadicuchiyama} Let $M$ be a non-positive sliced
  dyadic submartingale. Then $\nu=(\nu_I)_{I\in \mathcal{D}^4}$ with
  \[\nu_I=e^{M_I}\Delta^{\tmop{dy}} M_I\]
  is Carleson and there is the embedding from $H_{\tmop{dy}}^2$ to $\ell^2 (\nu)$: for all $f=u+iv \in H_{ \tmop{dy}}^2$ there holds
  \[ \sum_{I \in \mathcal{D}^4} \nu_I (u_I^2 + v_I^2) \leqslant   \| f \|^2_{H_{ \tmop{dy}}^2}. \]

\end{lemma}

 Our main result is below, Theorem \ref{corollaryuchiyamadyadic}. Note that we switched the signs, we prefer the formulation for non-negative supermartingale as opposed to non-positive submartingale (corresponding to the non-positive subharmonic function in the continuous case).

\begin{theorem}\label{corollaryuchiyamadyadic}
  \label{maintheoremuciyamastyle} Let $M$ be a bounded non-negative sliced
  dyadic supermartingale. Then setting $\mu=(\mu_I)_{I\in \mathcal{D}^4}$, $\frac{\mu_I}{|I|} = - \Delta^{\tmop{dy}}
  M_I$ there is the embedding from $H_{\tmop{dy}}^2$ to $\ell^2 (\mu)$: for all $f=u+iv \in H_{ \tmop{dy}}^2$ there holds
  \[ \sum_{I \in \mathcal{D}^4} \mu_I (u_I^2 + v_I^2) \leqslant e \| M
     \|_{\infty} \| f \|^2_{H_{ \tmop{dy}}^2} \]
  and the constant $e$ is sharp.
\end{theorem}

\paragraph{Relation to packing and testing conditions.}

 The classical complex embedding
theorem is stated in its geometric form, see \cite{C} .

\begin{theorem}[Carleson]
  \label{carlesonpacking}Let $\xi \in \mathbb{T}$ and $B (\xi, r)$ the ball in
  $\mathbb{C}$ with center $\xi$ and radius $r$. Define the Carleson box $C
  (\xi, r) = B (\xi, r) \cap \mathbb{D}$ and the Carleson intensity
  \[ I (\mu) = \sup \left\{ \frac{1}{r} \mu (C (\xi, r)) : r > 0, \xi \in
     \mathbb{T} \right\} . \]
  The embedding
  \[ \int_{\mathbb{D}} | f (z) |^2 \mathd \mu (z) \leqslant c (\mu) \| f
     \|_{H^2}^2 \]
  holds if and only if $I (\mu) < \infty$ with $c (\mu) \sim I (\mu) .$
\end{theorem}

There is an alternative formulation of Theorem \ref{corollaryuchiyama} that
resembles the Carleson Embedding theorem with packing condition, Theorem
\ref{carlesonpacking}. Let $S_z:\mathbb{D}\to\mathbb{D}$ denote the M\"{o}bius transform mapping 0 to $z\in \mathbb{D}$. Start with a subharmonic non-negative $\varphi$ and write
\[\varphi(z)= -\int_{\mathbb{C}}\log \frac1{|S_z(\xi)|}\Delta \varphi(\xi) \mathd A(\xi),\]
set  $\alpha (z)= -\Delta\varphi (z) \geqslant 0$ in $\mathbb{D}$ and
consider the Green's potential
\[ G (\alpha) (z) = \int_{\mathbb{D}} \log \frac{1}{| S_z (\xi) |} \alpha
   (\xi) \mathd A (\xi) . \]
Then $- \Delta G (\alpha) (z) = \alpha (z)$ and $\varphi (z) = - G
(\alpha) (z)$. So we may start with $\varphi $ or directly with the Green potential of a function $\alpha$. The assumption that
$\varphi (z)$ be bounded can now be rephrased as a packing condition on
$\alpha (z)$:
\[ G (\alpha) (z) = \int_{\mathbb{D}} \log \frac{1}{| S_z (\xi) |} \alpha
   (\xi) \mathd A (\xi) \leqslant 1. \]
Since $\log \frac{1}{| S_z (\xi) |}$ is the occupation time of the Brownian
motion in $\mathbb{D}$ begun at $z$ and stopped when hitting $\mathbb{T}$ \ so
that
\[ \int_{\mathbb{D}} \log \frac{1}{| S_z (\xi) |} \alpha (\xi) \mathd A (\xi)
   =\mathbb{E}^{\omega}_z \int^{\tau}_0 \alpha (B_s (\omega)) \mathd s, \]
this integral represents what of $\alpha (z)$ is picked up by Brownian
trajectories on average and it is therefore a packing condition in a direct analogy to dyadic packing conditions. In comparison, the dyadic
packing condition for the classical dyadic embedding reads $| I |^{- 1}
\sum_{J \in \mathcal{D} (I)} \mu_J \leqslant C \forall I$, which is again what
is picked up of the dyadic function $\mu_I$ when following dyadic trajectories
started at $I$. The reformulation of Theorem \ref{corollaryuchiyama} is thus:

\begin{lemma}[Uchiyama reformulated]
  Let $\alpha (z) \geqslant 0$ in $\mathbb{D}$ such that
  \[ \int_{\mathbb{D}} \log \frac{1}{| S_z (\xi) |} \alpha (\xi) \mathd A
     (\xi) \leqslant I (\alpha) \forall z \in \mathbb{D}. \]
  Then the measure
  \[ \mathd \mu (z) = \alpha (z) \log \frac{1}{| z |} \mathd A (z) \]
  produces the embedding $H^2 (\mathbb{D})$ in $L^2 (\mu)$ with constant $e$:
  \[ \int_{\mathbb{D}} | f (z) |^2 \mathd \mu (z) \leqslant e I (\alpha) \| f
     \|_{H^2}^2 . \]
\end{lemma}

\

Our main theorem can thus be stated resembling a dyadic
embedding, such as in Theorem \ref{maintheorempackingstyle} below. A 4-adic
sequence $(\mu_I)_{I \in \mathcal{D}^4 (I_0)}$ is balanced if
\begin{equation}
  \frac1{|I^x|}\sum_{J \in \mathcal{D}^4 (I^x)} \mu_J = \frac1{|I^y|}\sum_{J \in \mathcal{D}^4 (I^y)}
  \mu_J = \frac1{|I|}\sum_{J \in \mathcal{D}^4 (I)} \mu_J - \frac{\mu_I}{|I|} \label{condition-d=0}.
\end{equation}
This condition ensures that it stems from a sliced supermartingale $M$ via
\[M_I = \frac{1}{| I |} \sum_{J \in \mathcal{D}^4 (I)} \mu_J\]
so that $M_I-\frac{\mu_I}{|I|}=\frac12M_{I^{\star}_{+}}+\frac12M_{I^{\star}_{-}}$ for $\star=x,y$ separately.

One can thus restate our Theorem \ref{corollaryuchiyamadyadic} in this form below.

\begin{theorem}
  \label{maintheorempackingstyle}Let $\mu=(\mu_I)_{I\in \mathcal{D}^4}$ be a balanced non-negative sequence
  with the packing condition
  \[ \frac{1}{| I |} \sum_{J \in \mathcal{D}^4 (I)} \mu_J \leqslant C (\mu) \; \forall I\in \mathcal{D}^4,
  \]
  then there holds the embedding from $H_{\tmop{dy}}^2$ to $\ell^2 (\mu)$: for all $f=u+iv \in H_{\tmop{dy}}^2$ there holds
  \[ \sum_{I \in \mathcal{D}^4} \mu_I (u_I^2 + v_I^2) \leqslant e C (\mu) \| f
     \|^2_{H_{\tmop{dy}}^2} \]
  and the constant $e$ is sharp.
\end{theorem}

The condition (\ref{condition-d=0}) is required so that a sliced dyadic
supermartingale $M_I$ can be paired with the sequence $(\alpha_I)$. Our proof will
not work without this condition, see details in section \eqref{screwedup}.

\paragraph{Formulations using testing.}

From Theorem \ref{corollaryuchiyama}, one can then deduce via elementary
calculations that the testing form of the analytic Carleson embedding, a.k.a.
the reproducing kernel thesis for embedding theorems, holds. Such theorems test the statement on a representative class.

It is not hard to recover the testing condition of the
complex embedding once the Uchiyama lemma is proved. This was done in \cite{PTW}  via
the realization that an integral involving the reproducing kernels yields a non-positive
bounded subharmonic function. With a little more work, one can obtain very
good estimates, see \cite{PTW}:

\begin{theorem}[Petermichl, Treil, Wick]
  Let $\mu$ be any non-negative measure in $\mathbb{D}$ and for $\lambda \in \mathbb{D}$
  \[ \tilde{k}_{\lambda} (\xi) = \frac{(1 - | \lambda |^2)^{1 / 2}}{1 - \xi
     \bar{\lambda}} \]
  is the normalized reproducing kernel for the Hardy space $H^2 (\mathbb{D})$. If
  $\mu$ is Carleson in the testing sense with
  \begin{equation}
    A (\mu) = \sup_{\lambda \in \mathbb{D}} \| \tilde{k}_{\lambda} \|^2_{L^2 (\mu)} <
    \infty, \label{testing}
  \end{equation}
  then $H^2 (\mathbb{D})$ embeds in $L^2 (\mu)$ with constant $2 e :$
  \begin{equation}
    \int_{\mathbb{D}} | f (z) |^2 \mathd \mu (z) \leqslant 2 e A (\mu) \| f
    \|_{H^2}^2 . \label{embedding}
  \end{equation}
\end{theorem}

The condition (\ref{testing}) is clearly necessary for any estimate of the
type in (\ref{embedding}), as $\tilde{k}_{\lambda}$ are special cases for inequality
(\ref{embedding}), where $\| \tilde{k}_{\lambda} \|_{H^2}^2 = 1$ due to the
normalization. This type of condition is a so called testing condition as it
tests the estimate on a typical and representative enough class. The
multiplication of the true reproducing kernels of the Hardy space by the
bounded factor $(1 - | \lambda |^2)^{1 / 2}$ was essential to the argument in \cite{PTW}
and makes sense in the light of $L^2$ embeddings, even though these kernels
are not anymore `reproducing' in the classical sense.

\

In our case the dyadic reproducing kernel thesis is stated in the following theorem. Denote by $\tilde{k}^{\tmop{dy}}$ the normalized reproducing kernel of $H_{\tmop{dy}}^2$, that we derive below. Let again $\mu_I=-\Delta^{\tmop{dy}} M_I$ for $M$ a sliced non-negative supermartingale. We did not aim at good constants for this part of the estimate, but it might be hard (or impossible) to improve upon our estimate below.

\begin{theorem}\label{dyadictesting}
Assume that $\mu=(\mu_I)_{I\in \mathcal{D}^4}$ is balanced and \[\sup_I \|\tilde{k}^{\tmop{dy}}_I \|_{\ell^2_{\mu}}\le I(\mu).\] Then there is the embedding from  $H_{\tmop{dy}}^2$ to $\ell^2 (\mu)$: for all $f=u+iv \in H_{ \tmop{dy}}^2$ there holds \[ \sum_{I \in \mathcal{D}^4} \mu_I (u_I^2 + v_I^2) \leqslant 3e I(\mu)  \| f \|^2_{H_{ \tmop{dy}}^2}. \]

\end{theorem}

The Uchiyama lemma itself can be
interpreted as the analytic Carleson embedding theorem with Carleson packing
condition, which we explain below.

\subsection{A few words on our proofs}

\paragraph{Theorems \ref{corollaryuchiyamadyadic} and \ref{maintheorempackingstyle}.}
Both directions in our proof of Theorem \ref{corollaryuchiyamadyadic} or Theorem \ref{maintheorempackingstyle}, the bound itself and the optimality, are done by
a modification of the Bellman technique to adapt to analytic increments. An
extremal problem is set up that motivates the choice of a Bellman function.
Then the function we use is shown to have the required properties. To get the optimality we argue that certain aspects of this function cannot be improved. Even though the function itself is a supersolution to a Bellman equation, it is enough to get the optimality of the constant $e$.

\paragraph{Upper bound:}
Recall that a Bellman type function that gives the best estimate 4 in the
standard dyadic embedding is
\[ (F, f, M) \mapsto 4 F - \frac{4 f^2}{1 + M} \]
with its Hessian just barely negative semi-definite, a crucial feature for
proofs with Bellman functions. In its natural domain, where $0\leqslant M \leqslant 1$ and $F\geqslant f^2$, paired with the upper bound $4F$ and the non-negativity of the function and the derivative estimate in $M$ bounded below by $f^2$, this function proves in a by now standard way the constant 4 in the embedding theorem.

Our Bellman type function that gives the dyadic
analytic embedding is
\[ (F, u, v, M) \mapsto e F - e^{1 - M} (u^2 + v^2), \]
bearing obvious resemblance to the function used to prove Uchiyama's lemma. We
will see that it delivers the constant $e$. But it is easily checked that its
Hessian is indeed not negative semi-definite. In contrast, it has a restricted
concavity condition that respects dyadic analytic increments and this is what
allows us to have a constant $e$ better than 4. The estimate arises as a result of this restricted concavity paired with the non-negativity of the function, its upper estimate $eF$ and the estimate for the derivative in $M$ similar to above.

In the standard dyadic case, so-called mid point (of two points) concavity is implied by concavity as determined through its Hessian. This is in a sense a discrete analog of the comparison of an integral over a circle and the value in the center via the Laplacian. Notice that departing from two points and considering more points is not an obvious matter, see the insightful paper by Treil \cite{T}. In our case, there is a jump to four points, so the situation of the discrete and continuous case differs. Indeed, even though the Hessian matrix of our function shows the correct restricted concavity, it can fail for jumps to more than two points. As mentioned before, without the balancing condition on the sequence $\mu$ or the hypothesis that all processes are sliced, the required inequality that provides the dynamics fails.

\paragraph{Lower bound:}
For the optimality we modify the method of `reduction of a variable', an
unwritten argument from the late 90s {\cite{NTVcommunication}}, that recovered the
constant 4 as optimal in the dyadic embedding. The method consists of writing down an extremal problem for the estimate in the form of a Bellman function. By a so-called `reduction of a variable' the best constant for the original estimate enters the expression of the function. It is then shown that no function with the derived properties can exist with a constant better than 4. In our case it proves that $e$
can indeed not be improved, but again without giving an explicit counter example and via the set up of an extremal problem with a restricted concavity adapted to the dyadic analytic case. Due to the progression in the Bellman method, one might imagine different (lenghthier) ways to obtain this optimality, but we find our modification of this beautiful original argument by Nazarov, Treil and Volberg \cite{NTVcommunication} particularly minimalistic and elegant. Its extension to our case is interesting. We thank F. Nazarov for his help in recalling this argument in the standard dyadic case.


\paragraph{Lemma \ref{dyadicuchiyama}.}
To prove Lemma \ref{dyadicuchiyama} we use parts of the consideration of the proof for the upper bound of Theorem \ref{corollaryuchiyamadyadic}, but the argument itself appears to look different. It uses `convexity' instead of concavity' to deduce the estimate from a related Bellman function.

\paragraph{Theorem \ref{dyadictesting}.}
To deduce our Theorem \ref{dyadictesting}, we utilize Theorem \ref{corollaryuchiyamadyadic} directly instead of following the strategy in \cite{PTW}. Through the particular form of the reproducing kernel, we observe by direct calculation that the conditions of our Theorem \ref{corollaryuchiyamadyadic} are satisfied with a constant 3, if we assume the testing condition on the normalized reproducing kernels.

\section{Notation and formalism}\label{definitionsection}

Let $I_0 = [0, 1)$. As usual, the dyadic intervals are
\[ \mathcal{D} (I_0) = \{ 2^{- k} [0, 1) + n 2^{- k} : k, n \in \mathbb{N}_0,
   0 \leqslant n \leqslant 2^k - 1 \}. \]
   Let $I_-$ denote the left and $I_+$ the right half of the intervals $I$.
The dyadic Haar functions are $h_I = | I |^{- 1 / 2} (\chi_{I_+} - \chi_{I_-}) \;
\forall I \in \mathcal{D} (I_0)$ and form together with $ | I_0 |^{- 1 / 2} \chi_{I_0}$ an orthonormal basis of $L^2$. For any $L^2 (I_0)$ function $f$ we have
\[ f = \langle f \rangle_{I_0} + \sum_{I \in \mathcal{D} (I_0)} \langle f, h_I\rangle h_I  , \]
where $\langle \cdot \rangle_{I}$ denotes the average of a function over the interval $I$.
Notice that with $f_0=\langle f \rangle_{I_0}$,
\[ f_n = f_{0} + \sum_{\tmscript{\begin{array}{c}
     I \in \mathcal{D} (I_0)\\
     | I | > 2^{- n}
   \end{array}}} \langle f,h_I\rangle h_I \]
is measurable in the dyadic filtration $(\mathcal{F}_n)_{n \geqslant 0}$ where
$\mathcal{F}_n$ is the $\sigma$ algebra generated by $\mathcal{D}_n (I_0) = \{
I \in \mathcal{D} (I_0) : | I | = 2^{- n} \}$. $f=(f_n)_{n\geqslant 0}$ is a dyadic martingale
since $\mathbb{E} (f_{n + 1} - f_n | \mathcal{F}_n) = 0$. Its evaluation at time $n$ is $\sum_{I\in \mathcal{D}_n(I_0)} \langle f\rangle_I \chi_I$.

By the odd
respectively even dyadic intervals we mean
\[ \mathcal{D}_o (I_0) = \cup_{n \tmop{odd}} \mathcal{D}_n (I_0),
   \mathcal{D}_e (I_0) = \cup_{n \tmop{even}} \mathcal{D}_n (I_0) . \]
The 4-adic intervals are the even intervals
\[ \mathcal{D}^4 (I_0) = \mathcal{D}_e (I_0) = \{ 4^{- k} [0, 1) + n 4^{- k} : n, k \in
   \mathbb{N}_0, 0 \leqslant n \leqslant 4^k - 1 \} . \]

The martingales we consider in this paper are sliced dyadic martingales, so they can be written
as
\[ u_n = u_{0} + \sum_{\tmscript{\begin{array}{c}
     I \in \mathcal{D}_o (I_0)\\
     | I | > 2^{- n},
   \end{array}}} \langle u,h_I\rangle h_I = u_{0} + \sum_{\tmscript{\begin{array}{c}
     I \in \mathcal{D}_e (I_0)\\
     | I | = 2^{- n},
   \end{array}}} \langle u\rangle_I \chi_I\]
for even $n \geqslant 0$. Notice that in the notation with Haar functions, we have to use the odd intervals.
We use $u_I$ and $\langle u\rangle_I$ interchangeably where $u_I$ can be thought of $u$ sampled at the dyadic trajectory arriving in $I$ at time $\log|I|$.

If $I \in \mathcal{D}_e (I_0)$ then we identify $I_-
= I^y$ and $I_+ = I^x$. Notice that if $u$ and $v$ are sliced dyadic martingales then
\[ \frac{u_{I_+^x} + u_{I_-^x}}{2} = \frac{u_{I_+^y} + u_{I_-^y}}{2} = u,\qquad \frac{v_{I_+^x} +
   v_{I_-^x}}{2} =  \frac{v_{I_+^y} + v_{I_-^y}}{2} = v, \]
   so they are on the other hand 4-adic martingales with this additional restriction.

We call
\[ \Delta^{x} u_I = \frac{u_{I_+^x} - u_{I_-^x}}{2}, \Delta^y u_I = \frac{u_{I_+^y} - u_{I_-^y}}{2},
   \Delta^{x} v_I = \frac{v_{I_+^x} - v_{I_-^x}}{2}, \Delta^y v_I = \frac{v_{I_+^y} - v_{I_-^y}}{2}
\]
so that
\[ u_{I_\pm^x} = u_I \pm \Delta^x u_I, u_{I_\pm^y} = u_I \pm \Delta^y u_I, v_{I_\pm^x} = v_I
   \pm \Delta^x v_I, v_{I_\pm^y} = v_I \pm \Delta^y v_I. \]
We say that $v$ is the conjugate of $u$ if we have the Cauchy Riemann
equations
\[ \Delta^x u_I = \Delta^y v_I, \Delta^y u_I = - \Delta^x v_I. \]
We normalize the conjugate function so that $v_0=0$.

Notice that the operator from \cite{DP}
\[\mathcal{S}_0 : h_{I_{\pm}} \mapsto \pm h_{I_{\mp}}, I \in
\mathcal{D} (I_0) \backslash \{ I_0 \}, \mathcal{S}_0 : \chi_{I_0},h_{I_0} \mapsto 0\]
 delivers the conjugate: if $u$ is a sliced martingale then
$v =\mathcal{S}_0 (u)$ is its conjugate as
\begin{align*} \mathcal{S}_0(u)&=\mathcal{S}_0 \big(\langle u\rangle_{I_0}+ \sum_{I\in \mathcal{D}_o(I_0)}\langle u,h_I\rangle h_I\big)\\
&=\sum_{I\in \mathcal{D}_e(I_0)}\langle u,h_{I_-}\rangle \mathcal{S}_0h_{I_-}+\langle u,h_{I_+}\rangle \mathcal{S}_0h_{I_+}\\
&=\sum_{I\in \mathcal{D}_e(I_0)}-\langle u,h_{I_-}\rangle h_{I_+}+\langle u,h_{I_+}\rangle h_{I_-}\\
&=\sum_{I\in \mathcal{D}_e(I_0)}\langle v,h_{I_+}\rangle h_{I_+}+\langle u,h_{I_-}\rangle h_{I_-}=v.
\end{align*}
In the last line we used the Cauchy Riemann equations while recalling the convention $I^y=I_-$ and $I^x=I_+$ for even $I$, which imply $-\langle u,h_{I_-}\rangle =\langle v,h_{I_+}\rangle $ and $\langle u,h_{I_+}\rangle =\langle v,h_{I_-}\rangle $. Similarly $-\mathcal{S}_0(v)=u-u_0$.

A word of caution: the analytic `functions' should be viewed as a stochastic
process with dyadic analytic increments seen `from above', so that the
functions of the Hardy space in this setting can be identified with their $L^2$ closure
$f_{\infty} = u_{\infty} + i v_{\infty}$. One cannot just take any function in
$L^2 (I_0)$, extend it dyadically and try to find its conjugate - such a
function will not in general have a sliced dyadic representation. One may still think of $\mathbb{E} (u | \mathcal{F}_n)$ as a function of a variable
$t \in I_0$ via $\sum_{I \in \mathcal{D}_n} \langle u\rangle_I \chi_I$.

\

When working on $\mathbb{R}$ the dyadic intervals are, for example,
\[ \mathcal{D}  =\mathcal{D}(\mathbb{R})= \{ 2^{- k} [0, 1) + n 2^{- k} : k,n \in \mathbb{Z} \}, \]
the Haar series is over $n \in \mathbb{Z}$. The constant terms in the expressions and the moment conditions in the definition of $\mathcal{S}_0$ vanish.

\paragraph{Reproducing kernels.} As mentioned before, the dyadic Hardy space $H^2_{\tmop{dy}}$ consists of dyadic analytic $f = u + i v$ that are $L^2$ integrable. The inner product is inherited from $L^2$: $\langle f, g \rangle = \int f \bar{g}$. $H^2_{\tmop{dy}}$ is a reproducing kernel Hilbert space.

To find the reproducing kernel, let $f \in H_{\tmop{dy}}^2$. So for $I \in
\mathcal{D}^4$,
\[f_I = \langle f, \chi_I / | I | \rangle.\] 

Now recall that $f$ is sliced. Let $\mathcal{P}_o$ be the slicing operator,
projecting on the odd part: $\mathcal{P}_o \chi_{I_0}=\chi_{I_0}$,  $\mathcal{P}_o h_I = 0$ if $I \in \mathcal{D}_e$
 and $\mathcal{P}_o h_I = h_I $ if $I \in \mathcal{D}_o$. So
 \[\langle f,  \chi_I / | I | \rangle = \langle \mathcal{P}_o f,  \chi_I / | I | \rangle=
\langle f, \mathcal{P}_o  \chi_I / | I | \rangle .\]

The dyadic analytic projection is
\[ \mathcal{P}_{\tmop{dy}}^+ = \frac{1}{2} (\mathcal{I}+ i\mathcal{S}_0)
   \tmop{on} h_I, I \in \mathcal{D}_o (I_0), \mathcal{P}_{\tmop{dy}}^+
   =\mathcal{I} \tmop{on} \chi_{I_0} \]
with $\mathcal{P}_{\tmop{dy}}^+ f = f$ for $f \in H_{\tmop{dy}}^2$ since
\begin{align*}
\mathcal{P}_{\tmop{dy}}^+(u+iv)
&=\mathcal{P}_{\tmop{dy}}^+(u-u_0+iv)+u_0\\
&=u_0+\frac12(\mathcal{I}+ i\mathcal{S}_0)(u-u_0+iv)\\
&=u_0+\frac12(u-u_0+iv)+\frac{i}2(\mathcal{S}_0u+i\mathcal{S}_0v)\\
&=u+iv.
\end{align*}
The operator $\mathcal{P}_{\tmop{dy}}^+$ is self adjoint and
therefore
\[f_I = \langle \mathcal{P}_{\tmop{dy}}^+ f, \mathcal{P}_o\chi_I / | I | \rangle = \langle
f, \mathcal{P}_{\tmop{dy}}^+ \mathcal{P}_o\chi_I / | I | \rangle.\] 

Therefore
\[ k^{\tmop{dy}}_I = \mathcal{P}_{\tmop{dy}}^+ \mathcal{P}_o \chi_I / | I | \]
is the reproducing kernel. It can be computed explicitly: consider the dyadic
representation \ $\chi_I = \langle \chi_I \rangle_{I_0} + \sum_{J \in
\mathcal{D}: J \supsetneq I} \langle\chi_I, h_J\rangle h_J$. If $J= \hat{J}_{\pm}$ then denote $\sigma (J) = \pm$ and $J'$ the dyadic sibling of $J$ and $\hat{J}$ the dyadic parent.
Observe $\mathcal{S}_0 h_J = \sigma (J) h_{J'}$ for $J \in \mathcal{D}_o$.
\[ \mathcal{P}_{\tmop{dy}}^+\mathcal{P}_o \chi_I
= \langle \chi_I \rangle_{I_0} + \frac{1}{2}
\sum_{\tmscript{\begin{array}{c}
     J \in \mathcal{D}_o (I_0)\\
     J \supsetneq I,
   \end{array}}}
   \langle \chi_I, h_J\rangle h_J + \frac{i}{2}
   \sum_{\tmscript{\begin{array}{c}
     J \in \mathcal{D}_o (I_0)\\
     J \supsetneq I,
   \end{array}}}
   \sigma (J) \langle \chi_I, h_J\rangle h_{J'} . \]
Then we get
\begin{align*} k^{\tmop{dy}}_I &
=1
   + \frac{1}{2| I |}
   \sum_{\tmscript{\begin{array}{c}
     J \in \mathcal{D}_o (I_0)\\
     J \supsetneq I,
   \end{array}}}
   \langle {\chi_I}, h_J
   \rangle h_J + \frac{i}{2| I |}
   \sum_{\tmscript{\begin{array}{c}
     J \in \mathcal{D}_o (I_0)\\
     J \supsetneq I,
   \end{array}}}
   \sigma (J) \langle {\chi_I}, h_J \rangle h_{J'} .
\end{align*}
When working in $\mathbb{R}$ then
\begin{align*} k^{\mathbb{R},\tmop{dy}}_I
   &=  \frac{1}{2| I |}
   \sum_{\tmscript{\begin{array}{c}
     J \in \mathcal{D}_o (\mathbb{R})\\
     J \supsetneq I,
   \end{array}}}
   \langle {\chi_I}, h_J
   \rangle h_J + \frac{i}{2| I |}
   \sum_{\tmscript{\begin{array}{c}
     J \in \mathcal{D}_o (\mathbb{R})\\
     J \supsetneq I,
   \end{array}}}
   \sigma (J) \langle {\chi_I}, h_J \rangle h_{J'} .
\end{align*}

%
%
%
%
%


\

\section{The extremal problem.}

  To guide us in the proof of Theorems \ref{corollaryuchiyamadyadic} or \ref{maintheorempackingstyle}, we set up an extremal problem. It will on one hand
  guide our choice of a Bellman type function that we will use to prove the
  estimate. In addition, the extremal problem will be used directly to show
  the optimality.

  We set up the extremal problem as follows. In the domain
  \[ \Omega = \{ (\tmmathbf{F}, \tmmathbf{r}, \tmmathbf{i}, \tmmathbf{M}) :
     \tmmathbf{F} \geqslant \tmmathbf{r}^2 +\tmmathbf{i}^2, 0 \leqslant
     \tmmathbf{M} \leqslant 1 \} \subset \mathbb{R}^4 \]
  let
  \begin{align*} B (\tmmathbf{F}, \tmmathbf{r}, \tmmathbf{i}, \tmmathbf{M}) &= \sup_{M, f}
     \Big\{ \frac{1}{| I |} \sum_{J \in \mathcal{D}^4 (I)} \mu_J (u_J^2 +
     v_J^2)\Big| \\
     & u_I =\tmmathbf{r}, v_I =\tmmathbf{i}, M_I=\sum_{J \in \mathcal{D}^4
     (I)} \frac{\mu_J}{| I |} =\tmmathbf{M}, \langle u^2 + v^2 \rangle_I
     =\tmmathbf{F} \Big\},
  \end{align*}
  where the supremum runs on one hand over dyadic analytic $f = u + i v$ defined on $I\in \mathcal{D}^4$ with the values  $u_I =\tmmathbf{r}, v_I =\tmmathbf{i}$ and $\langle u^2 + v^2 \rangle_I=\tmmathbf{F} $ fixed. The latter is in the sense of the closure of the sliced dyadic martingales $u$ and $v$. 

  It runs on the other hand over sliced
  dyadic supermartingale $M$ stemming from a non-negative balanced 4-adic
  sequence $\mu_J : J \in \mathcal{D}^4(I)$ with
  \[ M_K=\sum_{J \in \mathcal{D}^4 (K)} \frac{\mu_J}{| K |}\leqslant 1\; \forall K\in \mathcal{D}^4  \] so that the value of the full sum $M_I=\tmmathbf{M}$ is fixed.
     $B$ as defined above depends upon $I$ but through scaling of the extremizer sequence, we can remove that dependence.

  \

  First, we show that the supremum in the definition of $B$ is not over the
  empty set, provided $(\tmmathbf{F}, \tmmathbf{r}, \tmmathbf{i},
  \tmmathbf{M}) \in \Omega$. For the variable $\tmmathbf{M}$ there is an
  obvious competitor, we can always choose the sequence $\mu_I =\tmmathbf{M} |
  I |, \mu_J = 0 \, \forall J \in \mathcal{D}^4 (I) \backslash \{ I \}$. For the
  other variables $\tmmathbf{F}, \tmmathbf{r}, \tmmathbf{i}$, we begin with the constant functions $\tilde{u}=\tmmathbf{r}$ and $\tilde{v}=\tmmathbf{i}$.
%
  From here, we choose the sliced orthogonal martingales
  below with constant $c$ to be determined:
  \[ u =\tilde{u} + c h_{I^x}^{\infty} + c h_{I^y}^{\infty},\qquad  v
     =\tilde{v} + c h_{I^x}^{\infty} - c h_{I^y}^{\infty}, \]
     where the Haar functions are normalized in $L^{\infty}$.
  Notice that $u_I =\tmmathbf{r}$ and $v_I =\tmmathbf{i}$. Now we compute
  \[ u^2 =\tilde{u}^2 + c^2 \chi_{I}  +
     2\tilde{u}c (h_{I^x}^{\infty} + h_{I^y}^{\infty}) \]
  \[ v^2 =\tilde{v}^2 + c^2 \chi_{I}  +
     2\tilde{v} c (-h_{I^x}^{\infty} + h_{I^y}^{\infty}) \]
  and thus since $\tilde{u}$ and $\tilde{v}$ are constant on $I$, we get $\langle u^2 + v^2 \rangle_I =\tmmathbf{r}^2 +\tmmathbf{i}^2 +
  2c^2 $. Since $\tmmathbf{F} \geqslant \tmmathbf{r}^2
  +\tmmathbf{i}^2$ we can choose $c$ so that $\tmmathbf{F}=\langle u^2 + v^2 \rangle_I$.
  Now we are ready to put the range
  \begin{equation}\label{rangeB} 0 \leqslant B (\tmmathbf{F}, \tmmathbf{r}, \tmmathbf{i}, \tmmathbf{M})
     \leqslant C\tmmathbf{F}. \end{equation}
  The lower estimate is because the supremum is of non-negative quantities and does not run over the empty set
  and the upper estimate is the belief that the theorem is true. The extremal problem is designed to make the most possible damage in every step along the way. Since we will use it to show optimality, it is important to observe that we start with a clean slate on each 4-adic interval $I$. That means, if analytic increments have been chosen between the intervals $I_0$ and a stopping time, then we are free to start over in any $I$, atom in the stopping filtration, and obtain a dyadic analytic continuation. So, the
  upper estimate in \eqref{rangeB} holds if and only if the theorem is true with constant $C.$
  We are looking for the best such $C$.

  We have a key inequality as follows. Let $0 \leqslant \tmmathbf{\mu}
  \leqslant \tmmathbf{M}$ then we claim the dynamics condition
  \begin{align}\label{fulldynamics}
    B (\tmmathbf{F}, \tmmathbf{r}, \tmmathbf{i}, \tmmathbf{M})  \geqslant &
    \tmmathbf{\mu} (\tmmathbf{r}^2 +\tmmathbf{i}^2)\\
      \nonumber& + \frac{1}{4} B (\tmmathbf{F}_-^x, \tmmathbf{r}_-^x,
    \tmmathbf{i}_-^x, \tmmathbf{M}_-^x) + \frac{1}{4} B (\tmmathbf{F}_+^x,
    \tmmathbf{r}_+^x, \tmmathbf{i}_+^x, \tmmathbf{M}_+^x)\\
      \nonumber& + \frac{1}{4} B (\tmmathbf{F}_-^y, \tmmathbf{r}_-^y,
    \tmmathbf{i}_-^y, \tmmathbf{M}_-^y) + \frac{1}{4} B (\tmmathbf{F}_+^y,
    \tmmathbf{r}_+^y, \tmmathbf{i}_+^y, \tmmathbf{M}_+^y)
  \end{align}
  as long as
  $ (\tmmathbf{F}, \tmmathbf{r}, \tmmathbf{i}, \tmmathbf{M})$,
     $(\tmmathbf{F}_-^x, \tmmathbf{r}_-^x, \tmmathbf{i}_-^x, \tmmathbf{M}_-^x)$,
     $(\tmmathbf{F}_+^x, \tmmathbf{r}_+^x, \tmmathbf{i}_+^x, \tmmathbf{M}_+^x)$,
     $(\tmmathbf{F}_-^y, \tmmathbf{r}_-^y, \tmmathbf{i}_-^y, \tmmathbf{M}_-^y)$,
     $(\tmmathbf{F}_-^y, \tmmathbf{r}_-^y, \tmmathbf{i}_-^y, \tmmathbf{M}_-^y)
     \in \Omega $
  and
  \[ \frac{1}{4} (\tmmathbf{F}_-^x +\tmmathbf{F}_+^x +\tmmathbf{F}_-^y
     +\tmmathbf{F}_+^y) =\tmmathbf{F}, \]
  \[ \frac{1}{2} (\tmmathbf{r}_-^x +\tmmathbf{r}_+^x) = \frac{1}{2}
     (\tmmathbf{r}_-^y +\tmmathbf{r}_+^y) =\tmmathbf{r} \tmop{with}
     \tmmathbf{r}_{\pm}^x =\tmmathbf{r} \pm \Delta^x \tmmathbf{r},
     \tmmathbf{r}_{\pm}^y =\tmmathbf{r} \pm \Delta^y \tmmathbf{r} \]
  \[ \frac{1}{2} (\tmmathbf{i}_-^x +\tmmathbf{i}_+^x) = \frac{1}{2}
     (\tmmathbf{i}_-^y +\tmmathbf{i}_+^y) =\tmmathbf{i} \tmop{with}
     \tmmathbf{i}_{\pm}^x =\tmmathbf{i} \mp \Delta^y \tmmathbf{r},
     \tmmathbf{i}_{\pm}^y =\tmmathbf{i} \pm \Delta^x \tmmathbf{r} \]
  \[ \frac{1}{2} (\tmmathbf{M}_-^x +\tmmathbf{M}_+^x) = \frac{1}{2}
     (\tmmathbf{M}_-^y +\tmmathbf{M}_+^y) =\tmmathbf{M}-\tmmathbf{\mu} \]
     with
     \[\tmmathbf{M}_{\pm}^x =\tmmathbf{M}-\tmmathbf{\mu} \pm \Delta^x \tmmathbf{M},
     \tmmathbf{M}_{\pm}^y =\tmmathbf{M}-\tmmathbf{\mu} \pm \Delta^y \tmmathbf{M}.\]

  Notice that we require Cauchy Riemann equations above on the increments of
  the pair $(\tmmathbf{r}, \tmmathbf{i})$. To ease our notation below, $\star \in \{ x, y \}, \sigma = \pm$. We show the dynamics condition  \eqref{fulldynamics}:

  \begin{align*}
    B (\tmmathbf{F}, \tmmathbf{r}, \tmmathbf{i}, \tmmathbf{M}) & \geqslant
    \sup \Big\{ \frac{1}{| I |} \sum_{J \in \mathcal{D}^4 (I)} \mu_J (\langle u
    \rangle_J^2 + \langle v \rangle_J^2)\Big|  \\
    &\quad \langle u
    \rangle_{I_{\sigma}^{\star}} =\tmmathbf{r}_{\sigma}^{\star}, \langle v
    \rangle_{I_{\sigma}^{\star}} =\tmmathbf{i}^{\star}_{\sigma},
     \langle u^2 + v^2 \rangle_{I_{\sigma}^{\star}}
    =\tmmathbf{F}_{\sigma}^{\star}, \sum_{J \in \mathcal{D}^4
    (I_{\sigma}^{\star})} \frac{\mu_J}{| I_{\sigma}^{\star} |}
    =\tmmathbf{M}_{\sigma}^{\star} \Big\}\\
    & =  \sup \Big\{ \sum_{\star, \sigma} \frac{1}{4 |I_{\sigma}^{\star}|} \sum_{J \in \mathcal{D}^4 (I^{\star}_{\sigma})} \mu_J (\langle u \rangle_J^2 + \langle v \rangle_J^2) +\tmmathbf{\mu} (\tmmathbf{r}^2
    +\tmmathbf{i}^2) \Big|\\
     &\quad  \langle u \rangle_{I_{\sigma}^{\star}}
    =\tmmathbf{r}_{\sigma}^{\star}, \langle v \rangle_{I_{\sigma}^{\star}}
    =\tmmathbf{i}^{\star}_{\sigma}, \langle u^2 + v^2
    \rangle_{I_{\sigma}^{\star}} =\tmmathbf{F}_{\sigma}^{\star},
    \sum_{J \in \mathcal{D}^4 (I_{\sigma}^{\star})} \frac{\mu_J}{|
    I_{\sigma}^{\star} |} =\tmmathbf{M}_{\sigma}^{\star} \Big\}\\
    & =  \frac{1}{4} \sum_{\star, \sigma} B (\tmmathbf{F}_{\sigma}^{\star},
    \tmmathbf{r}_{\sigma}^{\star}, \tmmathbf{i}_{\sigma}^{\star},
    \tmmathbf{M}_{\sigma}^{\star}) +\tmmathbf{\mu} (\tmmathbf{r}^2
    +\tmmathbf{i}^2) .
  \end{align*}
  We have used that the function $B$ itself does not depend upon $I$.

  Here are a few special cases of \eqref{fulldynamics}.
  If we choose $\tmmathbf{M}_{\sigma}^{\star} =\tmmathbf{M}-\tmmathbf{\mu}$
  and $\tmmathbf{F}_{\sigma}^{\star} =\tmmathbf{F},
  \tmmathbf{r}_{\sigma}^{\star} =\tmmathbf{r}, \tmmathbf{i}_{\sigma}^{\star}
  =\tmmathbf{i}$, then this becomes
  \begin{equation}\label{incrB} B (\tmmathbf{F}, \tmmathbf{r}, \tmmathbf{i}, \tmmathbf{M}) \geqslant B
     (\tmmathbf{F}, \tmmathbf{r}, \tmmathbf{i}, \tmmathbf{M}-\tmmathbf{\mu})
     +\tmmathbf{\mu} (\tmmathbf{r}^2 +\tmmathbf{i}^2) .
  \end{equation}
  If we choose $\tmmathbf{\mu}= 0$ 
%
%
  then
  \begin{equation}\label{concB} B (\tmmathbf{F}, \tmmathbf{r}, \tmmathbf{i}, \tmmathbf{M}) \geqslant
     \frac{1}{4} \sum_{\star, \sigma} B (\tmmathbf{F}_{\sigma}^{\star},
     \tmmathbf{r}_{\sigma}^{\star}, \tmmathbf{i}_{\sigma}^{\star},
     \tmmathbf{M}_{\sigma}^{\star}) . \end{equation}

  \section{Proof of Theorems \ref{corollaryuchiyamadyadic} and \ref{maintheorempackingstyle}}
  We first present the upper bounds, that is the embedding with constant $e$.

  \subsection{Upper Bound.}\label{upperbound} The idea is to find a function that has these
  properties $B$ has and we will see that this is enough to get the upper
  estimate. For our guess we modify the Bellman function used in the Uchiyama
  lemma $(\varphi, f) \mapsto e^{\varphi} | f |^2$ by switching to real
  variables and from the use of convexity to the use of concavity. Readers
  may already have noticed that we swapped the sign of the Laplacian as well.
  On the domain $\Omega = \{ \tmmathbf{r}^2 +\tmmathbf{i}^2 \leqslant
  \tmmathbf{F}, 0 \leqslant \tmmathbf{M} \leqslant 1 \}\subset \mathbb{R}^4$ we consider the
  function
  \[ \tilde{B} (\tmmathbf{F}, \tmmathbf{r}, \tmmathbf{i}, \tmmathbf{M}) =
     e\tmmathbf{F}- e^{1 -\tmmathbf{M}} (\tmmathbf{r}^2 +\tmmathbf{i}^2) . \]
  We claim that there holds for $(\tmmathbf{F}, \tmmathbf{r}, \tmmathbf{i},
  \tmmathbf{M}) \in \Omega$
  \begin{equation}
    0 \leqslant \tilde{B} (\tmmathbf{F}, \tmmathbf{r}, \tmmathbf{i},
    \tmmathbf{M}) \leqslant e\tmmathbf{F} \label{range}.
  \end{equation}
  If in addition  $(\tmmathbf{F}, \tmmathbf{r}, \tmmathbf{i},
  \tmmathbf{M}-\tmmathbf{\mu})\in \Omega$ then
  \begin{equation}
    \tilde{B} (\tmmathbf{F}, \tmmathbf{r}, \tmmathbf{i}, \tmmathbf{M}) -
    \tilde{B} (\tmmathbf{F}, \tmmathbf{r}, \tmmathbf{i},
    \tmmathbf{M}-\tmmathbf{\mu}) \geqslant \tmmathbf{\mu} (\tmmathbf{r}^2
    +\tmmathbf{i}^2) \label{derivative}
  \end{equation}
  and assuming that $(\tmmathbf{F}, \tmmathbf{r}, \tmmathbf{i}, \tmmathbf{M}),
  (\tmmathbf{F}_{\sigma}^{\star}, \tmmathbf{r}_{\sigma}^{\star},
  \tmmathbf{i}_{\sigma}^{\star}, \tmmathbf{M}_{\sigma}^{\star}) \in \Omega$
  with $\tmmathbf{\mu}=0$
  \[ \frac{\tmmathbf{F}_+^x +\tmmathbf{F}_-^x +\tmmathbf{F}_+^y
     +\tmmathbf{F}_-^y}{4} =\tmmathbf{F}, \frac{\tmmathbf{M}_+^{\star}
     +\tmmathbf{M}_-^{\star}}{2} =\tmmathbf{M}, \frac{\tmmathbf{r}_+^{\star}
     +\tmmathbf{r}_-^{\star}}{2} =\tmmathbf{r}, \frac{\tmmathbf{i}_+^{\star}
     +\tmmathbf{i}_-^{\star}}{2} =\tmmathbf{i} \]
  with Cauchy Riemann
  \[ \tmmathbf{r}_{\pm}^x =\tmmathbf{r} \pm \Delta^x \tmmathbf{r},
     \tmmathbf{r}_{\pm}^y =\tmmathbf{r} \pm \Delta^y \tmmathbf{r},
     \tmmathbf{i}_{\pm}^x =\tmmathbf{i} \mp \Delta^y \tmmathbf{r},
     \tmmathbf{i}_{\pm}^y =\tmmathbf{i} \pm \Delta^x \tmmathbf{r}, \]
then
  \begin{equation}
    \tilde{B} (\tmmathbf{F}, \tmmathbf{r}, \tmmathbf{i}, \tmmathbf{M})
    \geqslant \frac{1}{4} \sum_{\star, \sigma} \tilde{B}
    (\tmmathbf{F}_{\sigma}^{\star}, \tmmathbf{r}_{\sigma}^{\star},
    \tmmathbf{i}_{\sigma}^{\star}, \tmmathbf{M}_{\sigma}^{\star}).
    \label{concavity}
  \end{equation}
    To show (\ref{range}) we estimate
  \begin{align*} e\tmmathbf{F} & \geqslant e\tmmathbf{F}- e^{1 -\tmmathbf{M}}
     (\tmmathbf{r}^2 +\tmmathbf{i}^2)
     = \tilde{B}(\tmmathbf{F},\tmmathbf{r},\tmmathbf{i},\tmmathbf{M})\\&
     =e (\tmmathbf{F}- e^{-\tmmathbf{M}}
     (\tmmathbf{r}^2 +\tmmathbf{i}^2))
     \geqslant e (\tmmathbf{F}-
     (\tmmathbf{r}^2 +\tmmathbf{i}^2)) \geqslant 0.
  \end{align*}
  To show (\ref{derivative}) we calculate
  \[ \tilde{B} (\tmmathbf{F}, \tmmathbf{r}, \tmmathbf{i}, \tmmathbf{M}) -
     \tilde{B} (\tmmathbf{F}, \tmmathbf{r}, \tmmathbf{i},
     \tmmathbf{M}-\tmmathbf{\mu}) = e^{1 -\tmmathbf{M}} (\tmmathbf{r}^2
     +\tmmathbf{i}^2) (e^{\tmmathbf{\mu}} - 1) \geqslant \tmmathbf{\mu}
     (\tmmathbf{r}^2 +\tmmathbf{i}^2) . \]
  Inequality (\ref{concavity}) will now be shown by direct calculation.
  Disregarding the linear term of $\tilde{B}$ then dividing by $- e$ gives
  this inequality to prove:
  \begin{align*}
   e^{-\tmmathbf{M}} (\tmmathbf{r}^2 +\tmmathbf{i}^2)
   &\leqslant \frac{1}{4}
     e^{-\tmmathbf{M}_+^x} \big( {\tmmathbf{r}_+^x}^2
     {+\tmmathbf{i}_+^x}^2 \big)
     +\frac{1}{4} e^{-\tmmathbf{M}_-^x} \big(
     {\tmmathbf{r}_-^x}^2 {+\tmmathbf{i}_-^x}^2 \big) \\
     &\qquad +\frac{1}{4}
     e^{-\tmmathbf{M}_+^y} \big( {\tmmathbf{r}_+^y}^2 {+\tmmathbf{i}_+^y}^2
     \big)
     +\frac{1}{4} e^{-\tmmathbf{M}_-^y} \big( {\tmmathbf{r}_-^y}^2
     {+\tmmathbf{i}_-^y}^2 \big) .
     \end{align*}
  Rewriting the right hand side with increments gives this form:
  \begin{align*}
    e^{-\tmmathbf{M}} \big(\tmmathbf{r}^2 +\tmmathbf{i}^2\big)
    & \leqslant
    \frac{1}{4} e^{-\tmmathbf{M}_+^x} \big((\tmmathbf{r}+ \Delta^x
    \tmmathbf{r})^2 + (\tmmathbf{i}- \Delta^y \tmmathbf{r})^2\big) \\
    & \qquad + \frac{1}{4}
    e^{-\tmmathbf{M}_-^x} \big((\tmmathbf{r}- \Delta^x \tmmathbf{r})^2 +
    (\tmmathbf{i}+ \Delta^y \tmmathbf{r})^2\big)\\
    & \qquad   + \frac{1}{4} e^{-\tmmathbf{M}_+^y} \big((\tmmathbf{r}+ \Delta^y
    \tmmathbf{r})^2 + (\tmmathbf{i}+ \Delta^x \tmmathbf{r})^2\big) \\
    & \qquad  + \frac{1}{4}
    e^{-\tmmathbf{M}_-^y} \big((\tmmathbf{r}- \Delta^y \tmmathbf{r})^2 +
    (\tmmathbf{i}- \Delta^x \tmmathbf{r})^2\big) .
  \end{align*}
  This in turn is an inequality of quadratic forms acting on vectors
  $(\tmmathbf{r}, \tmmathbf{i}, \Delta^x \tmmathbf{r}, \Delta^y
  \tmmathbf{r})^t$. Let us bring the terms to the right hand side and sort out the terms.
  \begin{align*}
  & \Big( \frac{1}{4} e^{-\tmmathbf{M}_+^x} + \frac{1}{4}
     e^{-\tmmathbf{M}_-^x} + \frac{1}{4} e^{-\tmmathbf{M}_+^y} + \frac{1}{4}
     e^{-\tmmathbf{M}_-^y} - e^{-\tmmathbf{M}} \Big) (\tmmathbf{r}^2
     +\tmmathbf{i}^2), \\
  & \Big( \frac{1}{4} e^{-\tmmathbf{M}_+^x} + \frac{1}{4}
     e^{-\tmmathbf{M}_-^x} + \frac{1}{4} e^{-\tmmathbf{M}_+^y} + \frac{1}{4}
     e^{-\tmmathbf{M}_-^y} \Big) ((\Delta^x \tmmathbf{r})^2 + (\Delta^y
     \tmmathbf{r})^2), \\
  &  2 \Big( \frac{1}{4} e^{-\tmmathbf{M}_+^x} - \frac{1}{4}
     e^{-\tmmathbf{M}_-^x} \Big) \tmmathbf{r} \Delta^x \tmmathbf{r}, \\
  &  2 \Big( - \frac{1}{4} e^{-\tmmathbf{M}_+^x} + \frac{1}{4}
     e^{-\tmmathbf{M}_-^x} \Big) \tmmathbf{i} \Delta^y \tmmathbf{r}, \\
  &  2 \Big( \frac{1}{4} e^{-\tmmathbf{M}_+^y} - \frac{1}{4}
     e^{-\tmmathbf{M}_-^y} \Big) \tmmathbf{r} \Delta^y \tmmathbf{r}, \\
  &  2 \Big( \frac{1}{4} e^{-\tmmathbf{M}_+^y} - \frac{1}{4}
     e^{-\tmmathbf{M}_-^y} \Big) \tmmathbf{i} \Delta^x \tmmathbf{r}.
  \end{align*}
  The matrix
  \begin{equation} \label{mainmatrix}\Delta = \left( \begin{array}{cc} A&B\\ B^t&D\end{array} \right) \end{equation}
  with
  \begin{align*}
      A&=\left( \begin{array}{cc}
      \sum_{\star \sigma} e^{-\tmmathbf{M}_{\sigma}^{\star}} - 4
       e^{-\tmmathbf{M}}&0\\
       0&\sum_{\star \sigma} e^{-\tmmathbf{M}_{\sigma}^{\star}} - 4
       e^{-\tmmathbf{M}}
       \end{array} \right)\\
      B&=\left( \begin{array}{cc} e^{-\tmmathbf{M}_+^x} - e^{-\tmmathbf{M}_-^x}&e^{-\tmmathbf{M}_+^y} - e^{-\tmmathbf{M}_-^y}\\ e^{-\tmmathbf{M}_+^y} - e^{-\tmmathbf{M}_-^y}&
      -e^{-\tmmathbf{M}_+^x} + e^{-\tmmathbf{M}_-^x}\end{array} \right)\\
      C&=\left( \begin{array}{cc}
      \sum_{\star \sigma} e^{-\tmmathbf{M}_{\sigma}^{\star}}&0\\ 0&\sum_{\star \sigma} e^{-\tmmathbf{M}_{\sigma}^{\star}}
      \end{array} \right)
  \end{align*}
  acts on $(\tmmathbf{r}, \tmmathbf{i}, \Delta^x
  \tmmathbf{r}, \Delta^y \tmmathbf{r})^t$ as a bilinear form.
  We have to show it is positive semidefinite.
  Factoring out $e^{-\tmmathbf{M}}$
we get by naming
\[ \tmmathbf{M}_+^x -\tmmathbf{M}= d_1, \tmmathbf{M}_-^x -\tmmathbf{M}=
   - d_1, \tmmathbf{M}_+^y -\tmmathbf{M}=   d_2, \tmmathbf{M}_-^y
   -\tmmathbf{M}=  - d_2 \]
the matrix
\[ \frac{e^{-\tmmathbf{M}}}{4} \left( \begin{array}{cccc}
     - 4 + \Sigma & 0 & e^{ - d_1} - e^{ d_1} & e^{- d_2} - e^{
     d_2}\\
     0 & - 4 + \Sigma & e^{ - d_2} - e^{ d_2} & - e^{ - d_1} + e^{
     d_1}\\
     e^{ - d_1} - e^{ d_1} & e^{ - d_2} - e^{ d_2} & \Sigma & 0\\
     e^{ - d_2} - e^{ d_2} & - e^{ - d_1} + e^{ d_1} & 0 & \Sigma
   \end{array} \right) \]
where
\[ \Sigma = e^{ - d_1} + e^{ d_1} + e^{ - d_2} + e^{ d_2} . \]
To test if this matrix is positive semidefinite, we first observe that $- 4 + \Sigma \geqslant 0$ due to the convexity of the function $x\mapsto e^{-x}$, so the first two principal minors are non-negative. The third principal minor is
\[ \big( - 4 + \Sigma \big) \big[ \big( - 4 + \Sigma \big) \Sigma - (e^{ -
   d_1} - e^{ d_1})^2 - (e^{ - d_2} - e^{ d_2})^2 \big]. \]
This is after simplification
\[
  \big( - 4 + \Sigma \big) [- 4  (e^{- d_1} + e^{d_1}) - 4
  (e^{- d_2} + e^{d_2}) + 8 + 2 (e^{- d_1} + e^{d_1})
  (e^{- d_2} + e^{d_2})].
  \]
Recall that $- 4 + \Sigma \geqslant 0$ and observe $x\mapsto e^{-x}+e^x$ is bounded below by 2. So we consider the function $f(x_1,x_2)=-4x_1-4x_2+8+2x_1x_2=2(x_1-2)(x_2-2)$ for $x_1,x_2\geqslant 2$. It is elementary to show that the minimum is 0, attained on the boundary $f(x_1,2)=f(2,x_2)=0$, so that in turn the third principal minor is non-negative.
%
%
%

    It remains to compute the last principal minor, the determinant of the matrix itself, it is
    \[4\big(\frac{e^{-\tmmathbf{M}}}{4}\big)^4 (e^{\frac{d_1}{2}}-e^{-\frac{d_1}{2}})^4(e^{\frac{d_2}{2}}-e^{-\frac{d_2}{2}})^4\geqslant 0.\]
    So the matrix is positive semidefinite and (\ref{concavity}) is proved.

%


  \

  We are now ready to prove the desired estimate. Let $f_I = u_I + i v_I$ an
  analytic process in $H_{\tmop{dy}}^2$. For the first inequality, we invoke
  property (\ref{range}), followed by property (\ref{derivative}) for the
  second inequality and property (\ref{concavity}) for the third inequality.
  \begin{align*}
    &   | I_0 | e \langle u^2 + v^2 \rangle_{I_0}\\
    & \geqslant  | I_0 | \tilde{B} \Big( \langle u^2 + v^2 \rangle_{I_0},
    \langle u \rangle_{I_0}, \langle v \rangle_{I_0}, \frac{1}{| I_0 |}
    \sum_{J \in \mathcal{D}^4(I_0)} \mu_J \Big)\\
    & =  | I_0 | \tilde{B} \Big( \langle u^2 + v^2 \rangle_{I_0}, \langle u
    \rangle_{I_0}, \langle v \rangle_{I_0}, \frac{1}{| I_0 |} \sum_{J \in \mathcal{D}^4
    (I_0)} \mu_J \Big)\\
    & \qquad - | I_0 | \tilde{B} \Big( \langle u^2 + v^2
    \rangle_{I_0}, \langle u \rangle_{I_0}, \langle v \rangle_{I_0}, \sum_{J
    \in \mathcal{D}^4 (I_0)} \frac{\mu_J}{| I_0 |} - \frac{\mu_{I_0}}{| I_0 |} \Big)\\
    &  \qquad + | I_0 | \tilde{B} \Big( \langle u^2 + v^2 \rangle_{I_0}, \langle
    u \rangle_{I_0}, \langle v \rangle_{I_0}, \sum_{J \in \mathcal{D}^4 (I_0)}
    \frac{\mu_J}{| I_0 |} - \frac{\mu_{I_0}}{| I_0 |} \Big)\\
    & \geqslant  | I_0 | \frac{\mu_{I_0}}{| I_0 |} (\langle u
    \rangle_{I_0}^2 + \langle v \rangle_{I_0}^2) + | I_0 | \tilde{B} \Big(
    \langle u^2 + v^2 \rangle_{I_0}, \langle u \rangle_{I_0}, \langle v
    \rangle_{I_0}, \sum_{J \in \mathcal{D}^4 (I_0)} \frac{\mu_J}{| I_0 |} -
    \frac{\mu_{I_0}}{| I_0 |} \Big)\\
    & \geqslant   \mu_{I_0}(\langle u
    \rangle_{I_0}^2 + \langle v \rangle_{I_0}^2)\\
    &  \qquad +  \sum_{\tmscript{\begin{array}{cc}K \in \mathcal{D}^4 (I_0)\\ | K | = | I_0 | / 4 \end{array}}}  .  K | \tilde{B} \Big(
    \langle u^2 + v^2 \rangle_K, \langle u \rangle_K, \langle v \rangle_K,
    \frac{1}{| K |} \sum_{J \in \mathcal{D}^4 (K)} \mu_J \Big).
  \end{align*}

  Repeating this step and using both the upper and lower estimates in (\ref{range}) gives the desired upper bound. Observe that all 4 tuples
  above were always in the domain. For instance, we have used the restrictions that $M$ is
  sliced and that $f$ is dyadic analytic. Notice that we have the properties of the Bellman type
  function $\tilde{B}$ we used only for the values of such
  functions in the 4-adic tree.

  \

  \subsection{Lower bound.}
  We show that $e$ cannot be improved. We have seen in the previous subsection that for $f = u + i v \in H_{\tmop{dy}}^2$ and non-negative
  balanced sequences $\mu$ with $\frac{1}{| J |} \sum_{I \in \mathcal{D}^4 (J)}
  \mu_I \leqslant 1 \forall J$ there holds the upper bound
  \[ \sum_{I \in \mathcal{D}^4} \mu_I (u_I^2 + v_I^2) \leqslant e \| f
     \|^2_{H_{\tmop{dy}}^2} . \]
  We now show that the constant $e$ in the estimate cannot be improved by directly using the
  extremal problem we set up in the beginning. To do so, first consider we
  have (a better) estimate
  \[ \sum_{I \in \mathcal{D}^4} \mu_I (u_I^2 + v_I^2) \leqslant C \| f
     \|^2_{H_{ \tmop{dy}}^2} . \]
  Then we know that by construction $0 \leqslant B (\tmmathbf{F}, \tmmathbf{r}, \tmmathbf{i},
  \tmmathbf{M}) \leqslant C\tmmathbf{F}$, where $B$ is the function defined via the extremal problem. Let us now consider the reduction of a
  variable function
  \begin{equation} \label{defb}b (\tmmathbf{M},\tmmathbf{r}, \tmmathbf{i}) = \inf_{\tmmathbf{F}
     \geqslant \tmmathbf{r}^2 +\tmmathbf{i}^2} C\tmmathbf{F}- B (\tmmathbf{F},
     \tmmathbf{r}, \tmmathbf{i}, \tmmathbf{M})
     \end{equation}
  in the domain $\omega =\mathbb{R}^2 \times [0, 1]$. We will derive certain properties of $b$.

  Recall $B (\tmmathbf{F}, \tmmathbf{r},
  \tmmathbf{i}, \tmmathbf{M}) \leqslant C\tmmathbf{F}$ for all admissible constellations
  of $\tmmathbf{F}, \tmmathbf{r}, \tmmathbf{i},$ so in particular also for
  those that get close to the infimum in (\ref{defb}), so
  \begin{equation}
    b (\tmmathbf{M},\tmmathbf{r}, \tmmathbf{i}) \geqslant 0. \label{signb}
  \end{equation}
  The function $b$ also enjoys a bound from above:
  \begin{align*}
   0 &\leqslant B (\tmmathbf{F}, \tmmathbf{r}, \tmmathbf{i}, \tmmathbf{M}) =
     C\tmmathbf{F}- (C\tmmathbf{F}- B (\tmmathbf{F}, \tmmathbf{r},
     \tmmathbf{i}, \tmmathbf{M}))\\
      &\leqslant C\tmmathbf{F}- \inf_{\tmmathbf{F}
     \geqslant \tmmathbf{r}^2 +\tmmathbf{i}^2} (C\tmmathbf{F}- B
     (\tmmathbf{F}, \tmmathbf{r}, \tmmathbf{i}, \tmmathbf{M}))
     =
     C\tmmathbf{F}- b (\tmmathbf{M},\tmmathbf{r}, \tmmathbf{i}),
  \end{align*}
  so that for $(\tmmathbf{F}, \tmmathbf{r}, \tmmathbf{i},
  \tmmathbf{M}) \in \Omega$, we have $C\tmmathbf{F} \geqslant b (\tmmathbf{M},\tmmathbf{r}, \tmmathbf{i})$. In particular for $\tmmathbf{F}=\tmmathbf{r}^2+\tmmathbf{i}^2$ we obtain
  \begin{equation}
    C(\tmmathbf{r}^2+\tmmathbf{i}^2) \geqslant b (\tmmathbf{M},\tmmathbf{r}, \tmmathbf{i})
    \label{upperb} .
  \end{equation}
  The function $b$ inherits the `complex convexity' from the `complex concavity' of $B$: given
  $\tmmathbf{r}, \tmmathbf{i}, \tmmathbf{M}, \tmmathbf{r}_{\pm}^{\star},
  \tmmathbf{i}_{\pm}^{\star}, \tmmathbf{M}_{\pm}^{\star}$ so that
  ${\tmmathbf{r}}_{\pm}^{\star} =\tmmathbf{r} \pm \Delta^{\star} r,
  {\tmmathbf{i}}_{\pm}^{\star} =\tmmathbf{i} \pm \Delta^{\star} \tmmathbf{i}$
  with $\Delta^x \tmmathbf{r}= \Delta^y \tmmathbf{i}$ and $\Delta^y
  {\tmmathbf{r}}= - \Delta^x \tmmathbf{i}$ and $\tmmathbf{M}_{\pm}^{\star}
  ={\tmmathbf{M}} \pm \Delta^{\star} \tmmathbf{M}$. Choose
  $\tmmathbf{F}_{\sigma}^{\star}$ so that the infimum is almost attained in (\ref{defb}) in
  the definition of $b( \tmmathbf{M}_{\sigma}^{\star},\tmmathbf{r}_{\sigma}^{\star},
  \tmmathbf{i}_{\sigma}^{\star})$. For every $\varepsilon>0$ there is such a choice so that
  \begin{equation}\label{above} b ( \tmmathbf{M}_{\sigma}^{\star}\tmmathbf{r}_{\sigma}^{\star}, \tmmathbf{i}^{\star}_{\sigma}) \geqslant C\tmmathbf{F}_{\sigma}^{\star} - B
     (\tmmathbf{F}_{\sigma}^{\star}, \tmmathbf{r}_{\sigma}^{\star},
     \tmmathbf{i}_{\sigma}^{\star}, \tmmathbf{M}_{\sigma}^{\star}) - \varepsilon .
  \end{equation}
  Since the infimum runs over \ $\tmmathbf{F}^\star_\sigma \geqslant {\tmmathbf{r}^\star_\sigma} ^2
  +{\tmmathbf{i}^\star_\sigma}^2$, we know that $(\tmmathbf{F}_{\sigma}^{\star},
  \tmmathbf{r}_{\sigma}^{\star}, \tmmathbf{i}_{\sigma}^{\star},
  \tmmathbf{M}_{\sigma}^{\star}) \in \Omega$. Let us choose $\tmmathbf{F}=
  \frac{1}{4} (\tmmathbf{F}_+^x +\tmmathbf{F}_-^x +\tmmathbf{F}_+^y
  +\tmmathbf{F}_-^y) .$  By the convexity of $x\mapsto x^2$ we  have
  \[ \tmmathbf{r}^2
     \leqslant \frac{1}{2} ((\tmmathbf{r}_+^{\star})^2
     + (\tmmathbf{r}_-^{\star})^2), \qquad
   \tmmathbf{i}^2
     \leqslant \frac{1}{2} ((\tmmathbf{i}_+^{\star})^2
     + (\tmmathbf{i}_-^{\star})^2), \]
     so that using $\tmmathbf{F}_{\pm}^{\star} \geqslant
  (\tmmathbf{r}_{\pm}^{\star})^2 + (\tmmathbf{i}_{\pm}^{\star})^2$
  gives $\tmmathbf{F} \geqslant \tmmathbf{r}^2 +\tmmathbf{i}^2$ so that
  $(\tmmathbf{F}, \tmmathbf{r}, \tmmathbf{i}, \tmmathbf{M}) \in \Omega$.
  Now, using \eqref{concB} followed by \eqref{above}
  \begin{align*}
    b (\tmmathbf{M},\tmmathbf{r}, \tmmathbf{i})
    & \leqslant  C\tmmathbf{F}-
    B (\tmmathbf{F}, \tmmathbf{r}, \tmmathbf{i}, \tmmathbf{M})
     =  C \frac{1}{4} \sum_{\star, \sigma} \tmmathbf{F}_{\sigma}^{\star} -
    B (\tmmathbf{F}, \tmmathbf{r}, \tmmathbf{i}, \tmmathbf{M})\\
    & \leqslant  C \frac{1}{4} \sum_{\star, \sigma}
    \tmmathbf{F}_{\sigma}^{\star} - \frac{1}{4} \sum_{\star, \sigma} B
    (\tmmathbf{F}_{\sigma}^{\star}, \tmmathbf{r}_{\sigma}^{\star},
    \tmmathbf{i}_{\sigma}^{\star}, \tmmathbf{M}_{\sigma}^{\star})\\
    & \leqslant  \frac{1}{4} \sum_{\star, \sigma} b
    (\tmmathbf{M}^{\star}_{\sigma},\tmmathbf{r}_{\sigma}^{\star}, \tmmathbf{i}^{\star}_{\sigma}) + \varepsilon.
  \end{align*}
  Letting $\varepsilon \to 0$ we obtain
  \begin{equation}
    b (\tmmathbf{M},\tmmathbf{r}, \tmmathbf{i}) \leqslant \frac{1}{4}
    \sum_{\star, \sigma} b (\tmmathbf{M}^{\star}_{\sigma},\tmmathbf{r}_{\sigma}^{\star},
    \tmmathbf{i}^{\star}_{\sigma}), \label{convexityb}
  \end{equation}
  provided that all triples are in the domain $\omega$ and $\tmmathbf{r},
  \tmmathbf{i}, \tmmathbf{M}, \tmmathbf{r}_{\sigma}^{\star},
  \tmmathbf{i}_{\sigma}^{\star}, \tmmathbf{M}_{\sigma}^{\star}$ so that
  $\tmmathbf{r}_{\pm}^{\star} =\tmmathbf{r} \pm \Delta^{\star} r,
  \tmmathbf{i}_{\pm}^{\star} =\tmmathbf{i} \pm \Delta^{\star} \tmmathbf{i}$
  with $\Delta^x \tmmathbf{r}= \Delta^y \tmmathbf{i}$ and $\Delta^y
  \tmmathbf{r}= - \Delta^x \tmmathbf{i}$ and $\tmmathbf{M}_{\pm}^{\star}
  =\tmmathbf{M} \pm \Delta^{\star} \tmmathbf{M}$.

  There also holds in the same sense the increment property. Let
  $\tmmathbf{\mu} \in [0, \tmmathbf{M}]$ and $\tmmathbf{r}, \tmmathbf{i},
  \tmmathbf{M}$ given. For any $\varepsilon>0$ choose $\tmmathbf{F}$ so that $\tmmathbf{F} \geqslant
  \tmmathbf{r}^2 +\tmmathbf{i}^2$ and
  \[ b (\tmmathbf{M}-\tmmathbf{\mu},\tmmathbf{r}, \tmmathbf{i}) \geqslant
     C\tmmathbf{F}- B (\tmmathbf{F}, \tmmathbf{r}, \tmmathbf{i},
     \tmmathbf{M}-\tmmathbf{\mu}) - \varepsilon . \]
  Since $(\tmmathbf{F}, \tmmathbf{r}, \tmmathbf{i}, \tmmathbf{M}) \in \Omega$ and $\tmmathbf{\mu} \in [0, \tmmathbf{M}]$
  we may use \eqref{incrB} 
  to obtain
  \begin{align*} b (\tmmathbf{M},\tmmathbf{r}, \tmmathbf{i})
  &\leqslant C\tmmathbf{F}- B
     (\tmmathbf{F}, \tmmathbf{r}, \tmmathbf{i}, \tmmathbf{M})
   \leqslant  C\tmmathbf{F}- B (\tmmathbf{F}, \tmmathbf{r}, \tmmathbf{i},
     \tmmathbf{M}-\tmmathbf{\mu}) -\tmmathbf{\mu} (\tmmathbf{r}^2
     +\tmmathbf{i}^2) \\
 &\leqslant b (\tmmathbf{M}-\tmmathbf{\mu},\tmmathbf{r}, \tmmathbf{i}) -\tmmathbf{\mu} (\tmmathbf{r}^2
     +\tmmathbf{i}^2) + \varepsilon .
  \end{align*}
  Letting $\varepsilon \rightarrow 0$ gives
  \begin{equation}
    b (\tmmathbf{M},\tmmathbf{r}, \tmmathbf{i}) \leqslant b (\tmmathbf{M}-\tmmathbf{\mu},\tmmathbf{r},
    \tmmathbf{i}) -\tmmathbf{\mu} (\tmmathbf{r}^2
    +\tmmathbf{i}^2), \label{derivativeb}
  \end{equation}
  provided both triples are in the domain $\omega$.

  In addition, the Bellman function $B$ has the scaling $B (t^2 \tmmathbf{F},
  t\tmmathbf{r}, t\tmmathbf{i}, \tmmathbf{M}) = t^2 B (\tmmathbf{F},
  \tmmathbf{r}, \tmmathbf{i}, \tmmathbf{M})$ and this is inherited by $b$:
  \begin{align*}
   b (\tmmathbf{M},t\tmmathbf{r}, t\tmmathbf{i})
   &= \inf_{t^2 \tmmathbf{F}
     \geqslant (t\tmmathbf{r})^2 + (t\tmmathbf{i})^2} C t^2 \tmmathbf{F}- B
     (t^2 \tmmathbf{F}, t\tmmathbf{r}, t\tmmathbf{i}, \tmmathbf{M}) \\
     &= t^2
     \inf_{\tmmathbf{F} \geqslant \tmmathbf{r}^2 +\tmmathbf{i}^2}
     C\tmmathbf{F}- B (\tmmathbf{F}, \tmmathbf{r}, \tmmathbf{i}, \tmmathbf{M})
     = t^2 b (\tmmathbf{M},\tmmathbf{r}, \tmmathbf{i})   .
  \end{align*}
  From this homogeneity of degree 2 it follows that $b$ can be defined on
  $[0, 1] \times  B(0,1)$ and then extended. We now argue that $b (\tmmathbf{M},1, 0) = b (\tmmathbf{M},\tmmathbf{r}, \tmmathbf{i})$ for all
  $(\tmmathbf{r}, \tmmathbf{i}) \in B(0,1)$. Indeed, obtain
  $(\tmmathbf{r}, \tmmathbf{i}) \in B(0,1)$ through a rotation of the point $(1, 0)\in B(0,1)$
  around the origin through the angle $\theta \in [0, 2 \pi)$ so that
  $\Omicron_{\theta} (1, 0) = (\tmmathbf{r}, \tmmathbf{i})$. Let $f = u + i v$
  be a dyadic analytic map and the sequence $(\mu_I)$ for $\tmmathbf{M}$ so that
  the supremum in the definition of $B$ is almost attained:
  \[ B (\tmmathbf{F}, 1, 0, \tmmathbf{M}) \leqslant \frac{1}{| I |} \sum_{J
     \in \mathcal{D}^4 (I)} \mu_J (\langle u \rangle_J^2 + \langle v
     \rangle_J^2) + \varepsilon \]
  where $\langle u \rangle_I = 1, \langle v \rangle_I = 0, \langle u^2 + v^2
  \rangle_I =\tmmathbf{F}, \frac{1}{| I |} \sum_{J \in \mathcal{D}^4 (I)} \mu_J
  =\tmmathbf{M}$. Any rotation of the values of $f$ will not affect the
  variables $\tmmathbf{F}, \tmmathbf{M}$ or the value of the sum. The rotation of $f$ will
  retain dyadic analyticity because the orthogonality of the increments is preserved under a rotation. It follows that the rotated $f$ together with the original $\mu$ will almost attain the supremum in the definition of $B (\tmmathbf{F}, \tmmathbf{r}, \tmmathbf{i}, \tmmathbf{M})$.
  Thus $b$ is constant on $\{ \tmmathbf{M} \} \times B(0,1)  $ for
  $0 \leqslant \tmmathbf{M} \leqslant 1$.

  \

  In this light by scaling our function has the form
  \[ b (\tmmathbf{M},\tmmathbf{r}, \tmmathbf{i}) = (\tmmathbf{r}^2
     +\tmmathbf{i}^2) \Phi (\tmmathbf{M}), \]
  with $\Phi (\tmmathbf{M}) = b (1, 0, \tmmathbf{M})$. We thus obtain for our
  previous properties the following:
  The size properties (\ref{signb}) and (\ref{upperb}) become
  \begin{equation}\label{sizePhi}
  0\leqslant\Phi(\tmmathbf{M})\leqslant C.
  \end{equation}

  Property (\ref{derivativeb}) becomes $(\tmmathbf{r}^2 +\tmmathbf{i}^2) \Phi
  (\tmmathbf{M}) \leqslant (\tmmathbf{r}^2 +\tmmathbf{i}^2) \Phi
  (\tmmathbf{M}-\tmmathbf{\mu}) -\tmmathbf{\mu} (\tmmathbf{r}^2
  +\tmmathbf{i}^2)$ and thus
  \begin{equation} \label{derivativephi}\Phi (\tmmathbf{M}) \leqslant \Phi (\tmmathbf{M}-\tmmathbf{\mu})
     -\tmmathbf{\mu}. \end{equation}

  Notice that $B (\tmmathbf{F}, \tmmathbf{r}, \tmmathbf{i},0) = 0$ because
  when $\tmmathbf{M}= 0$ at a certain time then it remains 0 for all future
  times. This implies that $b ( \tmmathbf{r}, \tmmathbf{i},0) =
  \inf_{\tmmathbf{F} \geqslant \tmmathbf{r}^2 +\tmmathbf{i}^2} C\tmmathbf{F}=
  C (\tmmathbf{r}^2 +\tmmathbf{i}^2)$ so that $\Phi (0) = C \neq 0$. Since $\Phi$ is
  strictly decreasing and $\Phi\geqslant 0$ on $[0,1]$, we have $\Phi > 0$ on $(0,1)$.

  As for the convexity property (\ref{convexityb}) we get
  \begin{equation}\label{convexityphi}(\tmmathbf{r}^2 +\tmmathbf{i}^2) \Phi (\tmmathbf{M}) \leqslant
     \frac{1}{4} \sum_{\star, \sigma} ((\tmmathbf{r}_{\sigma}^{\star})^2 +
     (\tmmathbf{i}^{\star}_{\sigma})^2) \Phi (\tmmathbf{M}^{\star}_{\sigma}) . \end{equation}
  The function $\Phi$ might not be smooth. Mollifying in $\tmmathbf{M}$ in a standard way with a bell shaped $\varphi_{\varepsilon}$ supported on $[-\varepsilon,\varepsilon]$ gives us a smooth function $\tilde{\Phi}_{\varepsilon}(\tmmathbf{M})$ defined on $[\varepsilon,1-\varepsilon]$. Notice that the so obtained function has $\tilde{\Phi}_{\varepsilon}(\tmmathbf{M})>0$ on $[\varepsilon,1-\varepsilon]$. We deduce from (\ref{derivativephi}) and (\ref{convexityphi}) the same inequalities for
  $$\tilde{b}_{\varepsilon}(\tmmathbf{M},\tmmathbf{r},\tmmathbf{i})=(\tmmathbf{r}^2+\tmmathbf{i}^2)\tilde{\Phi}_{\varepsilon}(\tmmathbf{M})$$ in $[\varepsilon,1-\varepsilon] \times \mathbb{R}^2 $. The properties (\ref{sizePhi}) and (\ref{derivativephi}) imply for $\tmmathbf{M} \in [\varepsilon,1-\varepsilon]$
  $$0 < \tilde{\Phi}_{\varepsilon}(\tmmathbf{M})< C, \qquad \tilde{\Phi}_{\varepsilon}'(\tmmathbf{M})\leqslant -1.$$

   We develop $\tilde{b}_{\varepsilon}(\tmmathbf{M},\tmmathbf{r},\tmmathbf{i})$
   by Taylor. Let $X =
  (\tmmathbf{M},\tmmathbf{r}, \tmmathbf{i})^t$ and $X_{\pm}^{\star} = X +
  \Delta_{\pm}^{\star} X = X \pm \Delta^{\star} X$ with $\Delta^x
  \tmmathbf{r}= \Delta^y \tmmathbf{i}$ and $\Delta^y \tmmathbf{r}= - \Delta^x
  \tmmathbf{i}$. For some $0 \leqslant \xi^\star_\sigma \leqslant 1$ we have
  \begin{equation}\label{discreteconvexity} \sum_{\star, \sigma} \tilde b_{\varepsilon} (X_{\sigma}^{\star}) = 4 \tilde b_{\varepsilon} (X) + 0 + \frac{1}{2}
     \sum_{\star, \sigma} \langle \mathd^2 \tilde b_{\varepsilon} (X + \xi^\star_\sigma \Delta_{\sigma}^{\star}
     X) \Delta_{\sigma}^{\star} X, \Delta_{\sigma}^{\star} X \rangle .
  \end{equation}
  Observing Cauchy-Riemann, we have for $\Delta_{\sigma}^{\star} X$ the following vectors
  \begin{align*}
  \Delta_{+}^{x} X&=\left(\begin{array}{c}
     \Delta^x \tmmathbf{M}\\
      \Delta^x \tmmathbf{r}\\
      - \Delta^y \tmmathbf{r}
    \end{array} \right) , \,
    \Delta_{-}^{x} X=\left(\begin{array}{c}
     -\Delta^x \tmmathbf{M}\\
      -\Delta^x \tmmathbf{r}\\
       \Delta^y \tmmathbf{r}
    \end{array} \right),\\
  \Delta_{+}^{y} X&=\left(\begin{array}{c}
     \Delta^y \tmmathbf{M}\\
      \Delta^y \tmmathbf{r}\\
       \Delta^x \tmmathbf{r}
    \end{array} \right),\,
    \Delta_{-}^{y} X=\left(\begin{array}{c}
     -\Delta^y \tmmathbf{M}\\
      -\Delta^y \tmmathbf{r}\\
       -\Delta^x \tmmathbf{r}
    \end{array} \right).
  \end{align*}

  Writing the very right hand side of (\ref{discreteconvexity}) as a quadratic form on the increments $(\Delta^x \tmmathbf{M},\Delta^y \tmmathbf{M},\Delta^x \tmmathbf{r},\Delta^y \tmmathbf{r})^t$ and then
  letting all the increments go to 0, this converges to a non-negative definite quadratic form acting on  $(\mathd^x \tmmathbf{M},\mathd^y \tmmathbf{M},\mathd^x \tmmathbf{r},\mathd^y \tmmathbf{r})^t$
  with matrix

  \[\left( \begin{array}{cccc}
       \frac{\partial^2\tilde b_{\varepsilon}}{\partial\tmmathbf{M}^2} (X) &
       0 &
       \frac{\partial^2\tilde b_{\varepsilon}}{\partial\tmmathbf{M}\partial\tmmathbf{r}}(X) &
       - \frac{\partial\tilde b_{\varepsilon}} {\partial\tmmathbf{M}\partial\tmmathbf{i}}(X)\\
       0 &
       \frac{\partial^2\tilde b_{\varepsilon}}{\partial\tmmathbf{M}^2} (X) &
       \frac{\partial^2\tilde b_{\varepsilon}}{\partial\tmmathbf{M}\partial\tmmathbf{i}}(X) &
       \frac{\partial^2\tilde b_{\varepsilon}}{\partial\tmmathbf{M}\partial\tmmathbf{r}} (X)\\
       \frac{\partial^2\tilde b_{\varepsilon}}{\partial\tmmathbf{M}\partial\tmmathbf{r}} (X) &
       \frac{\partial^2\tilde b_{\varepsilon}}{\partial\tmmathbf{M}\partial\tmmathbf{i}} (X) &
       \frac{\partial^2\tilde b_{\varepsilon}}{\partial\tmmathbf{r}^2} (X) +
       \frac{\partial^2\tilde b_{\varepsilon}}{\partial\tmmathbf{i}^2} (X) &
       0\\
       - \frac{\partial^2\tilde b_{\varepsilon}}{\partial\tmmathbf{M}\partial\tmmathbf{i}} (X) &
       \frac{\partial^2\tilde b_{\varepsilon}}{\partial\tmmathbf{M}\partial\tmmathbf{r}}(X) &
       0 &
       \frac{\partial^2\tilde b_{\varepsilon}}{\partial\tmmathbf{r}^2} (X) +
       \frac{\partial^2\tilde b_{\varepsilon}}{\partial\tmmathbf{i}^2} (X)
     \end{array} \right). \]

  Computing these derivatives for the function $(\tmmathbf{r}^2
  +\tmmathbf{i}^2) \tilde \Phi_{\varepsilon} (\tmmathbf{M})$ gives the following continuous analog of (\ref{convexityphi}):
  \[ \left( \begin{array}{cccc}
       (\tmmathbf{r}^2 +\tmmathbf{i}^2) \tilde \Phi_{\varepsilon}'' (\tmmathbf{M}) & 0 &
       2\tmmathbf{r} \tilde \Phi_{\varepsilon}' (\tmmathbf{M}) & - 2\tmmathbf{i} \tilde \Phi_{\varepsilon}'
       (\tmmathbf{M})\\
       0 & (\tmmathbf{r}^2 +\tmmathbf{i}^2) \tilde \Phi_{\varepsilon}'' (\tmmathbf{M}) &
       2\tmmathbf{i} \tilde \Phi_{\varepsilon}' (\tmmathbf{M}) & 2\tmmathbf{r} \tilde \Phi_{\varepsilon}'
       (\tmmathbf{M})\\
       2\tmmathbf{r} \tilde \Phi_{\varepsilon}' (\tmmathbf{M}) & 2\tmmathbf{i} \tilde \Phi_{\varepsilon}' (\tmmathbf{M})
       & 4\tilde \Phi_{\varepsilon} (\tmmathbf{M}) & 0\\
       - 2\tmmathbf{i} \tilde \Phi_{\varepsilon}' (\tmmathbf{M}) & 2\tmmathbf{r} \tilde \Phi_{\varepsilon}'
       (\tmmathbf{M}) & 0 & 4 \tilde \Phi_{\varepsilon} (\tmmathbf{M})
     \end{array} \right) \geqslant 0. \]
  We compute the four principal minors $M_1 - M_4$, starting in the right
  lower corner.
  \begin{eqnarray*}
    M_1 & = & 4 \tilde \Phi_{\varepsilon} (\tmmathbf{M})\\
    M_2 & = & 16 \tilde \Phi_{\varepsilon} (\tmmathbf{M})^2\\
    M_3 & = & 16 (\tmmathbf{r}^2 +\tmmathbf{i}^2) \tilde \Phi_{\varepsilon} (\tmmathbf{M})
    [\tilde \Phi_{\varepsilon}'' (\tmmathbf{M}) \tilde \Phi_{\varepsilon} (\tmmathbf{M}) - \tilde \Phi_{\varepsilon}' (\tmmathbf{M})^2]\\
    M_4 & = & 16 (\tmmathbf{r}^2 +\tmmathbf{i}^2)^2 [\tilde \Phi_{\varepsilon}'' (\tmmathbf{M})
    \tilde \Phi_{\varepsilon} (\tmmathbf{M}) - \tilde \Phi_{\varepsilon}' (\tmmathbf{M})^2]^2
  \end{eqnarray*}
  For the matrix to be positive semidefinite, we get the requirements
  \[ \tilde \Phi_{\varepsilon} (\tmmathbf{M}) \geqslant 0 \infixand \tilde \Phi_{\varepsilon}'' (\tmmathbf{M}) \tilde \Phi_{\varepsilon}
     (\tmmathbf{M}) - \tilde \Phi_{\varepsilon}' (\tmmathbf{M})^2 \geqslant 0. \]
  Recalling that $\tilde\Phi_{\varepsilon}(1-\varepsilon)>0$ and thus $\tilde\Phi_{\varepsilon}(\tmmathbf{M})>0$ on $[\varepsilon,1-\varepsilon]$, the latter is equivalent to the convexity of
  \[ \Upsilon_{\varepsilon} : \tmmathbf{M} \mapsto \log \tilde \Phi_{\varepsilon} (\tmmathbf{M}). \]
  Indeed
  \[ \frac{\partial^2 \Upsilon_{\varepsilon}  (\tmmathbf{M})} { \partial \tmmathbf{M}^2}=
  \frac{\partial^2 \log \tilde \Phi_{\varepsilon} (\tmmathbf{M})} { \partial \tmmathbf{M}^2}
     = - \frac{1}{\tilde \Phi_{\varepsilon} (\tmmathbf{M})^2} \tilde \Phi_{\varepsilon}' (\tmmathbf{M})^2 +
     \frac{1}{\Phi_{\varepsilon} (\tmmathbf{M})} \tilde \Phi_{\varepsilon}'' (\tmmathbf{M}) \]
   and multiplying by $\tilde \Phi_{\varepsilon} (\tmmathbf{M})^2$ gives $0 \leqslant \tilde\Phi_{\varepsilon}''
  (\tmmathbf{M}) \tilde\Phi_{\varepsilon} (\tmmathbf{M}) - \tilde\Phi_{\varepsilon}' (\tmmathbf{M})^2$, if and only if
  $\Upsilon_{\varepsilon} $ is convex. $\tilde\Phi_{\varepsilon}' (\tmmathbf{M}) \leqslant - 1$ so
  \[ \Upsilon_{\varepsilon}' (\tmmathbf{M}) = \frac{1}{\tilde\Phi_{\varepsilon} (\tmmathbf{M})} \tilde\Phi_{\varepsilon}'
     (\tmmathbf{M}) \leqslant - \frac{1}{\tilde\Phi_{\varepsilon} (\tmmathbf{M})} = -
     \frac{1}{{e^{\Upsilon_{\varepsilon} (\tmmathbf{M})}} } . \]
  Let us observe
  \[a_{\varepsilon} := \Upsilon_{\varepsilon} (1-\varepsilon)\leqslant \log C, \qquad\Upsilon_{\varepsilon}' (1-\varepsilon) \leqslant - \frac{1}{e^{a_{\varepsilon}}}.\] When $s  <  1-\varepsilon$ we get by
  convexity
  \[\Upsilon_{\varepsilon} (s) \geqslant \Upsilon_{\varepsilon} (1-\varepsilon) + \Upsilon_{\varepsilon}' (1-\varepsilon) (s - 1+\varepsilon) \geqslant a_{\varepsilon} + \frac{1}{e^{a_{\varepsilon}}} (1 - s-\varepsilon)\]
   and so
   \[\Upsilon_{\varepsilon} (\varepsilon) \geqslant a_{\varepsilon} + \frac{1}{e^{a_{\varepsilon}}}(1-2\varepsilon) \geqslant 1+2\varepsilon a_{\varepsilon}-2\varepsilon,\] which in turn implies
   \[\tilde\Phi_{\varepsilon} (\varepsilon) \geqslant e \cdot e^{2\varepsilon(a_{\varepsilon}-1)}.\]
   Thanks to \eqref{derivativephi} applied to $\tmmathbf{M}=1$ and $\varepsilon\leqslant \tmmathbf{\mu} \leqslant 2\varepsilon$ we know that $\Phi(\tmmathbf{M})\geqslant \varepsilon$ for $1-2\varepsilon \leqslant \tmmathbf{M}\leqslant \varepsilon$. Since we mollify with a bell shaped, symmetric, $L^1$ normalized function, we get $\tilde\Phi_{\varepsilon}(1-\varepsilon)\geqslant \frac12 \varepsilon$.
   Therefore, since $a_{\varepsilon}=\log \tilde\Phi_{\varepsilon}(1-\varepsilon)$, we get $|a_{\varepsilon}|\leqslant|\log\frac12\varepsilon |$ and thus $\varepsilon a_{\varepsilon}\to 0$ as $\varepsilon \to 0$.
   Letting $\varepsilon \to 0$  gives  $e^{2\varepsilon(a_{\varepsilon}-1)}\to 1$. Observing that $\Phi(0)\geqslant \tilde\Phi_{\varepsilon} (\varepsilon)$ gives the constant $e$ is optimal.

  \section{Dyadic Uchiyama}

To prove Lemma \ref{dyadicuchiyama} we use the properties of the function $\tilde{B}$ to obtain similar properties for the function
\[
(\tmmathbf{r},\tmmathbf{i},\tmmathbf{M})\mapsto e^{\tmmathbf{M}}(\tmmathbf{r}^2
  +\tmmathbf{i}^2).
\]
 Let either $M_I$ be a non-positive sliced
  dyadic submartingale and set $\mu_I=|I|\Delta^{\tmop{dy}}M_I$ or start with balanced $\mu_I
\geqslant 0$ for $I \in \mathcal{D}^4 (I_0)$ and put
\[ M_I = - \frac{1}{| I |} \sum_{J \in \mathcal{D}^4 (I)} \mu_J, \]
so that $M_I$ non-positive with $\Delta^{\tmop{dy}} M_I = \frac{\mu_I}{| I |}.$
Compute the full dyadic Laplacian
\begin{align*}
&\Delta_{\tmop{full}}^{\tmop{dy}} e^{M_I} (u_I^2 + v_I^2)\\
    := & \sum_{\star, \sigma} \frac{1}{4} e^{M_{I_{\sigma}^{\star}}}
  (u_{I_{\sigma}^{\star}}^2 + v_{I_{\sigma}^{\star}}^2) - e^{M_I} (u_I^2 +
  v_I^2)\\
   = & \sum_{\star, \sigma} \frac{1}{4} e^{M_{I_{\sigma}^{\star}}}
  (u_{I_{\sigma}^{\star}}^2 + v_{I_{\sigma}^{\star}}^2) - e^{M_I} (u_I^2 +
  v_I^2) + e^{M_I + \frac{\mu_I}{| I |}} (u_I^2 + v_I^2) - e^{M_I +
  \frac{\mu_I}{| I |}} (u_I^2 + v_I^2)\\
   \geqslant & \sum_{\star, \sigma} \frac{1}{4} e^{M_{I_{\sigma}^{\star}}}
  (u_{I_{\sigma}^{\star}}^2 + v_{I_{\sigma}^{\star}}^2) - e^{M_I +
  \frac{\mu_I}{| I |}} (u_I^2 + v_I^2) + \frac{\mu_I}{| I |} e^{M_I} (u_I^2 +
  v_I^2)\\
   \geqslant & \frac{\mu_I}{| I |} e^{M_I} (u_I^2 + v_I^2).
\end{align*}
In this estimate we used properties of $\tilde{B}$, namely \eqref{derivative} and \eqref{concavity}.
Now by the above, applied to $I=I_0$, we obtain 
 \begin{align*}
      0
     \leqslant & | I_0 | e^{M_{I_0}} \langle u^2 + v^2 \rangle_{I_0}\\[.5em]
     = & | I_0 | e^{M_{I_0}} \langle u^2 + v^2 \rangle_{I_0} +
     \mu_{I_0} e^{M_{I_0}} (u_{I_0}^2 + v_{I_0}^2) - \mu_{I_0} e^{M_{I_0}}
     (u_{I_0}^2 + v_{I_0}^2)\\[.5em]
     \leqslant & \sum_{\star, \sigma} |{{I_0}_{\sigma}^{\star}}| e^{M_{{I_0}_{\sigma}^{\star}}}
     (u_{{I_0}_{\sigma}^{\star}}^2 + v_{{I_0}_{\sigma}^{\star}}^2) - \mu_{I_0}
     e^{M_{I_0}} (u_{I_0}^2 + v_{I_0}^2).
   \end{align*}
Iterating this estimate, we obtain
   \[ 0\leqslant  \sum_{\tmscript{\begin{array}{cc}J \in \mathcal{D}^4 (I_0)\\ | J | = 2^{- n - 2} \end{array}}} | J | e^{M_J} (u_J^2 + v_J^2) - \sum_{\tmscript{\begin{array}{cc}J \in \mathcal{D}^4 (I_0)\\ | J | > 2^{- n - 2} \end{array}}} \mu_J e^{M_J} (u_J^2 + v_J^2) .\]
Letting $n \rightarrow \infty$ and observing $1 \geqslant e^{M_J}$ we obtain
the estimate
\[ \int_{I_0} | f |^2 \geqslant \sum_{J \in \mathcal{D}^4 (I_0)} \mu_J e^{M_J}
   (u_J^2 + v_J^2) \]
so that there is an embedding with the measure $\mu_J e^{M_J}$.

\section{The reproducing kernel of the dyadic Hardy space}

For Theorem \ref{dyadictesting} it is more convenient to work in $\mathbb{R}$.
The testing condition is
\[ \sum_{K \in \mathcal{D}^4(\mathbb{R}) } | \langle {k}^{\mathbb{R},\tmop{dy}}_I \rangle_{K} |^2 \mu_K \leqslant \|{k}^{\mathbb{R},\tmop{dy}}_I\|_{H^2_{\tmop{dy}}}^2, \]
where we now compute the norm on the right in order to normalize the reproducing kernel. Recall that

\[ {k}^{\mathbb{R},\tmop{dy}}_I =\frac1{|I|} \mathcal{P}_{\tmop{dy}}^+\mathcal{P}_o \chi_I \]
with
\[ {k}^{\mathbb{R},\tmop{dy}}_I =
    \frac{1}{2} \sum_{J \in \mathcal{D}_o(\mathbb{R}) : J
   \supsetneq I} \langle\frac{\chi_I}{|I|}, h_J\rangle h_J + \frac{i}{2} \sum_{J \in \mathcal{D}_o(\mathbb{R}) : J
   \supsetneq I} \sigma (J) \langle\frac{\chi_I}{|I|}, h_J\rangle h_{J'}.  \]
Then the reproducing kernel property is
\[
f(K) := \langle f \rangle_K
=
\big\langle
f,
k^{\mathbb{R},\tmop{dy}}_K\
\big\rangle
\]
for all $f \in H^2_{\tmop{dy}}$ and 4-adic $K \subseteq \mathbb{R}$. 
In order to compute the norm of the kernel, observe that for $J \in \mathcal{D}_o(\mathbb{R}), J\supsetneq I$ we have $h_{J'}(I)=0 $.
Using the reproducing kernel property we get
\[
\|{k}^{\mathbb{R},\tmop{dy}}_I\|_{H^2_{\tmop{dy}}}^2
=
{k}^{\mathbb{R},\tmop{dy}}_I(I)
=
\frac12 \sum_{J \in \mathcal{D}_o(\mathbb{R}), J\supsetneq I}
\frac1{|J|}
=
\frac12 \sum_{k\geqslant 0}
\frac1{2\cdot 4^k |I|}
=\frac1{3|I|}.
\]
Therefore set \[ \tilde{k}^{\mathbb{R},\tmop{dy}}_I=\sqrt{3|I|}{k}^{\mathbb{R},\tmop{dy}}_I,\]
   so that the testing condition becomes
   \begin{equation}\label{normalizedcondition} \sum_{K \in \mathcal{D}^4(\mathbb{R}) } | \tilde{k}^{\mathbb{R},\tmop{dy}}_I (K) |^2 \mu_K \leqslant 1. \end{equation}
We now compute $| \tilde{k}^{\mathbb{R},\tmop{dy}}_I (K) |^2$ for any 4-adic $K\subset I$. In addition note that, for $J \in \mathcal{D}_o(\mathbb{R})$ and $J\supsetneq I$ we have $h_{J'}(K)=0 $ and $h_{J}(K)=h_{J}(I)$ for all 4-adic $K \subseteq I$. It follows again by the reproducing kernel property

\[\tilde{k}^{\mathbb{R},\tmop{dy}}_I(K)
=
\tilde{k}^{\mathbb{R},\tmop{dy}}_I(I)
=
\langle \tilde{k}^{\mathbb{R},\tmop{dy}}_I, k^{\mathbb{R},\tmop{dy}}_I\rangle
=
\frac
{\langle k^{\mathbb{R},\tmop{dy}}_I, k^{\mathbb{R},\tmop{dy}}_I\rangle}
{\|k^{\mathbb{R},\tmop{dy}}_I\|}
=
\|{k}^{\mathbb{R},\tmop{dy}}_I\|_{H^2_{\tmop{dy}}}
=
\frac{1}{\sqrt{3|I|}}.
\]
Hence, (\ref{normalizedcondition}) immediately implies
\[\frac1{|I|}\sum_{K\in \mathcal{D}^4(I)}\mu_K\leqslant 3,\]
for which we then get an embedding with constant $3e$ for the testing formulation.

\section{Slicing or balancing condition}\label{screwedup}
In this section we remark that slicing cannot be removed for the upper estimate. Indeed, removing the slicing condition on the variable $\tmmathbf{M}$ alone is not possible for the Bellman function we used. If we remove the balancing condition on the sequence $\mu_J$ or the slicing on $M$, we obtain, similar to subsection \ref{upperbound} a quadratic form that has to be positive definite in order for the restricted concavity condition to hold.
Dividing each entry by $e^{-\tmmathbf{M}}$ of the matrix $\Delta$ from \eqref{mainmatrix} acting on $(\tmmathbf{r}, \tmmathbf{i}, \Delta^x \tmmathbf{r},
\Delta^y \tmmathbf{r})^t$, we get by naming
\[ \tmmathbf{M}_+^x -\tmmathbf{M}= d + d_1, \tmmathbf{M}_-^x -\tmmathbf{M}= d
   - d_1, \tmmathbf{M}_+^y -\tmmathbf{M}= - d + d_2, \tmmathbf{M}_-^y
   -\tmmathbf{M}= - d - d_2 \]
and
   \[ \Sigma = e^{- d - d_1} + e^{- d + d_1} + e^{d - d_2} + e^{d + d_2} = e^{- d}
   (e^{- d_1} + e^{d_1}) + e^d (e^{- d_2} + e^{d_2})  \]
the matrix
\[\frac{e^{-\tmmathbf{M}}}{4}
\left( \begin{array}{cc}
D&E\\E^t&F
\end{array}\right)
\]
with
\begin{align*}
    D&=\left(\begin{array}{cc}
-4+\Sigma&0\\0&-4+\Sigma
\end{array}\right)\\
    E&=\left(\begin{array}{cc}
    e^{- d}(e^{ - d_1} - e^{  d_1}) & e^d(e^{ - d_2} - e^{
     d_2})\\
     e^d(e^{ - d_2} - e^{ d_2}) & - e^{-d}(e^{ - d_1} + e^{ +
     d_1})
    \end{array}\right)\\
    F&=\left(\begin{array}{cc}
\Sigma&0\\0&\Sigma
\end{array}\right).
\end{align*}

We will observe that the third minor (modulo the constant factor) can be negative definite. It
is
\begin{equation} \label{thirdminorall}\big( - 4 + \Sigma \big) \left[ \big( - 4 + \Sigma \big) \Sigma - (e^{- d -
   d_1} - e^{- d + d_1})^2 - (e^{d - d_2} - e^{d + d_2})^2 \right],
   \end{equation}
which expands to
\begin{align*}
G(d,d_1,d_2)&=\big(-4+e^{- d} (e^{- d_1} + e^{d_1}) + e^d
  (e^{- d_2} + e^{d_2})\big)F(d,d_1,d_2)\\
F(d,d_1,d_2)&=- 4 e^{- d} (e^{- d_1} + e^{d_1}) - 4 e^d
  (e^{- d_2}+ e^{d_2})\\
  &\qquad  + 4 e^{- 2 d} + 4 e^{2 d} + 2 (e^{- d_1} + e^{d_1})
  (e^{- d_2} + e^{d_2}).
\end{align*}
 If we swap
the roles of $d_1$ and $d_2$ we see that we may assume without loss of
generality that $d \geqslant 0$. Likewise $d_1, d_2 \geqslant 0$ without loss
of generality. We have as domain
\[ 0 \leqslant \tmmathbf{M}, \tmmathbf{M}+ d + d_1, \tmmathbf{M}+ d - d_1,
   \tmmathbf{M}- d + d_2, \tmmathbf{M}- d - d_2 \leqslant 1, \]
so that we get for $\delta_{\tmmathbf{M}} = \tmop{dist} (\tmmathbf{M},
\partial I_0)$ that $d + d_1, d + d_2 \leqslant \delta_{\tmmathbf{M}}
\leqslant 1 / 2$. Thus $0 \leqslant d \leqslant \delta_{\tmmathbf{M}}, d_{1,
2} \leqslant \delta_{\tmmathbf{M}} - d$.
By elementary calculation we find that for some $d$ near 0 and $d_1=0$, $d_2=\frac12-d$ the expression
\eqref{thirdminorall} can be negative.

\

\begin{center}
\includegraphics[height=6cm]{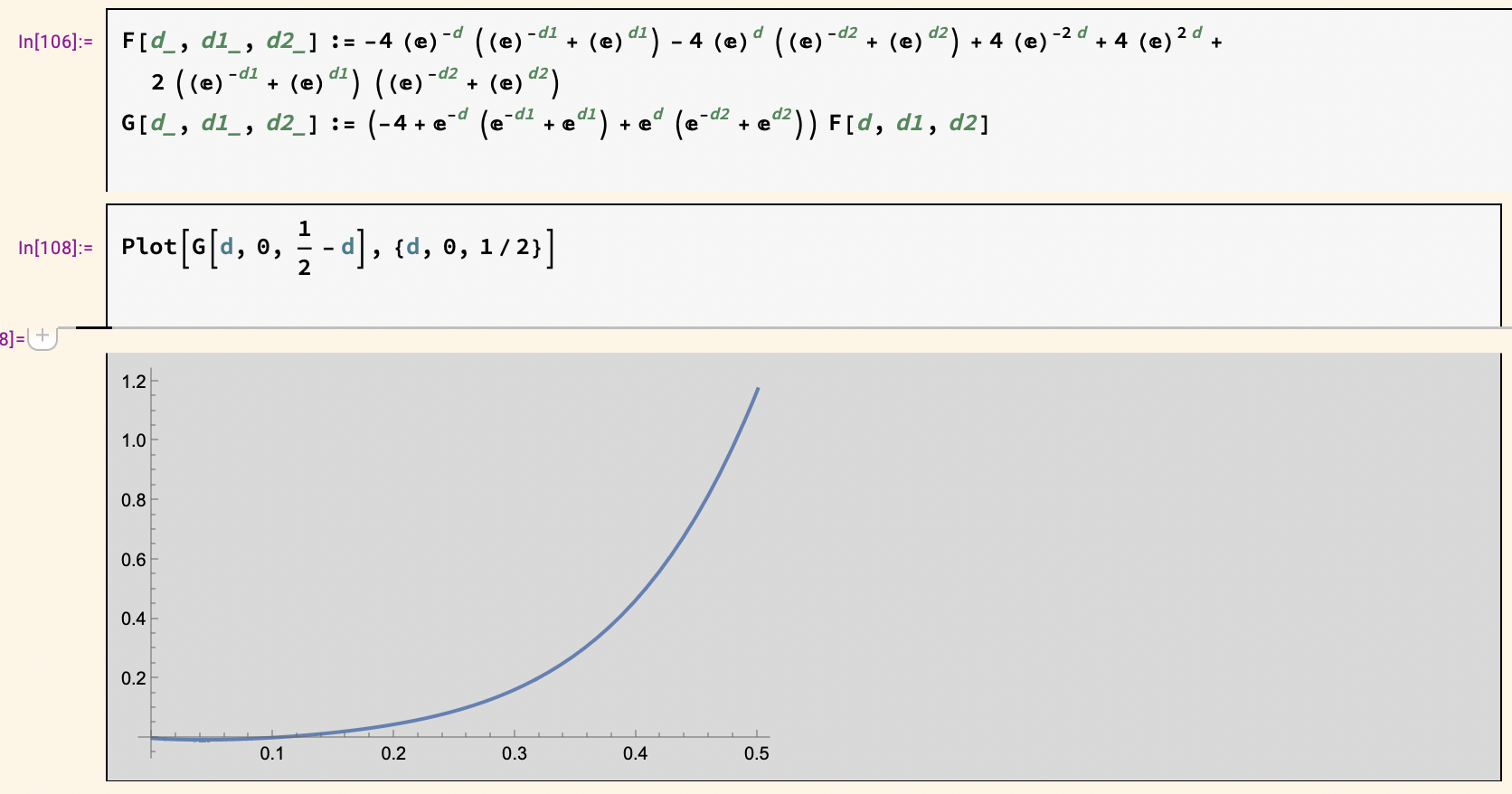}
\end{center}

\

\

\

\

\

\

\end{document}